\documentclass[10pt]{article}
\usepackage{anysize}
\marginsize{2.0cm}{2.0cm}{2.0 cm}{2.0 cm}
\usepackage{setspace}
\usepackage{graphicx}
\usepackage{color}
\usepackage{epstopdf}
\usepackage[titletoc]{appendix}
\usepackage{graphicx}
\usepackage{latexsym,amsfonts,amsmath,amssymb}
\newtheorem{remark}{Remark}

\newtheorem{condition}{Condition}
\newtheorem{theorem}{Theorem}
\newtheorem{lemma}{Lemma}

\newtheorem{definition}{Definition}
\newtheorem{proposition}{Proposition}
\newtheorem{example}{Example}

\def\RR{\mathbb{R}}

\def\PP{\mathsf{P}}

\def\II{\mathbb{I}}

\def\EE{\mathsf{E}}
\def\PP{\mathsf{P}}

\begin{document}

\title{Primal-dual programs for the constrained optimal impulse control: discounted model}


\author{Alexey Piunovskiy\thanks{Corresponding author. Department of Mathematical Sciences, University of
Liverpool, Liverpool, U.K. E-mail: piunov@liv.ac.uk.}
\and Yi
Zhang \thanks{School of Mathematics, University of Birmingham, Edgbaston,
Birmingham,
B15 2TT, U.K. Email: y.zhang.29@bham.ac.uk} }

\maketitle

\par\noindent\textbf{Abstract}. This paper studies constrained optimal impulse control problems of a deterministic system described by a (semi)flow, where the performance measures are the discounted total costs including both the costs incurred with  applying  impulses as well as running costs.  We formulate the relaxed problem and the associated primal convex programs in measures together with the dual programs, and establish the relevant duality results. As an application, we formulate and justify a general procedure for obtaining optimal $(J+1)$-mixed strategies for the original impulse control problems. This procedure is illustrated with a solved example.

\bigskip


\section{Introduction}\label{sec1}
This paper studies optimal impulse control problems of a deterministic system described by a (semi)flow, where the performance measures are the discounted total costs including both the lump sum cost associated with   applying impulses as well as running costs. The decision maker aims at minimizing one performance measure, subject to constraints on the other similar performance measures given by different impulse cost functions and running cost rates. 

We devote this paper to Professor A. Plakhov on the occasion of his 65th birthday. Let  us start with comments on this paper in relation to the remarkable contributions made by Professor A. Plakhov  in this area, while postponing the general literature review until a later part of this section. In absence of constraints and discounting, together with Professor A. Plakhov, we investigated this problem in \cite{PiunovskiySasha:2018}, where the dynamic programming approach was studied.  The first step of investigation in that paper was to reformulate the impulse control problem as a Markov decision process (MDP), leading to the so called optimality equation in integral form, and after that, optimality equation in differential form (close to the so called quasi-variational inequalities) was established, and the connection between these two equations was studied in detail. This work admitted natural extensions. One is for more general system dynamics, and this was carried out together with Professor A. Plakhov in \cite{dpp} for optimal impulse control problems of piecewise deterministic processes. Another possibility is to consider problems with (functional) constraints. For handling such problems, the convex analytic approach, also called the linear programming method, is quite convenient.  This method is based on reformulating the original problem as a static optimization problem in suitably defined occupation measures. This optimization problem can be relaxed to a convex program in measures after one imposes the constraint given by the characteristic equation satisfied by occupation measures. A fundamental issue in this convex analytic approach was to investigate the equivalence of this convex program and the original impulse control problem, known as the absence of relaxation gap. This was the main result in  \cite{pz1,PiunovskiyZhangJMAA:2021b}, where we acknowledged the initial input from Professor A. Plakhov. We mention that corresponding to the two types of optimality equations for unconstrained problems, one may deal with two types of occupation measures, leading to two convex program formulations, which we may again call of integral form and of differential form. The first one was investigated in \cite{pz1}, and the second one was investigated in  \cite{PiunovskiyZhangJMAA:2021b}.  After the absence of relaxation gap is showed, the next natural step is to formulate the dual program of the convex program, and investigate the absence of duality gap. This was done in \cite{pz2}, but for unconstrained optimal impulse control problem with total undiscounted costs.  

The aim of the present paper is to extend and complement the results in \cite{pz2} to problems with constraints and discounted costs, for which we further investigate the convex analytic approach by focusing on the duality issues and formulate procedures for finding optimal strategies in specific forms, the so called  $(J+1)$-mixed strategies with $J$ being the number of constraints. Note that, without further conditions, the relaxed problem of the original impulse control problem is in space of infinite measures, even in presence of discounting. This is in contrast with situations encountered in convex analytic approach for (classic) optimal gradual control problems. The primal convex program is in the occupation measures similar to those considered in \cite{pz1}. One may formulate its dual program, and the crude but otherwise useful observation is the absence of the duality gap, which follows from the relevant results in \cite{b42}. The dual functional is in the space of non-negative Lagrange multipliers. It is desirable (and is required in our general procedure to solve constrained optimal impulse problems) to compute the value of the dual functional for each Lagrange multiplier by solving an unconstrained impulse control problem with a new cost function modified with the Lagrange multipliers.  To guarantee that these two values coincide, we impose suitable conditions on the cost functions, which also imply that the relaxed problem is equivalent to a (so called second) convex program in finite measures. Having said the above, our main contributions lie in the suitable formulation of the dual program of the second convex program, proving the absence of the duality gap, showing that the dual functional of the relaxed problem can be computed by solving an unconstrained problem  using the dynamic programming method, and presenting with justifications the general procedure for computing the optimal  $(J+1)$-mixed optimal strategies for the impulse control problems. We also demonstrate the procedure by solving in closed form an example that can be understood in the context of fluid queueing systems. We note that while the existence of optimal mixed strategies over $J+1$ deterministic stationary strategies, which was not addressed in the aforementioned works, follows from the recent paper \cite{sicon24}, here we showed that the mixture could be done over deterministic stationary strategies which are uniformly optimal for the unconstrained problem modulated by  optimal Lagrange  multipliers.

Finally, we briefly review the other relevant and recent literature on impulse control problems.
 In the absence of constraints, the standard method of investigation of optimal impulse control problems is dynamic programming, leading to the so called quasi-variational inequalities \cite{bar,chris,dpp}, see e.g., \cite{bar,chris,dpp,jelito,PiunovskiySasha:2018,ppt}. The Pontryagin maximum principle and similar methods were developed in \cite{aru,Miller:2003}. Another approach that is especially convenient in handling problems with constraints is the convex analytic approach. This method is based on reformulating the original problem as linear or convex programs in the space of suitably defined occupation measures. Such measures and linear programs were introduced and investigated in \cite{cho,helmes,pz1,pz2} in different settings. The convex analytic approach can be also developed for classical (gradual) optimal control problems, see \cite{gai}.  
A little more details  about the cited works are as follows.\begin{itemize}
\item The state space was a subset $\RR^n$ and the dynamical system was defined through an ordinary differential equation in \cite{aru,bar,gai,helmes,Miller:2003,ppt}. The more general (e.g., metric) state space appeared in \cite{cho,chris,dpp,jelito,PiunovskiySasha:2018,pz1,pz2}. The system was stochastic in \cite{cho,chris,helmes,jelito}.
\item The most common objective to be minimized was the total (expected) cost in \cite{dpp,PiunovskiySasha:2018,ppt,pz1,pz2}, and the special case of discounted cost was also studied in \cite{bar,chris,helmes,jelito,PiunovskiySasha:2018}. The (expected) terminal cost and the long-run average cost were the objectives in \cite{aru,cho,Miller:2003} and \cite{gai} correspondingly.
\item Functional constraints appeared in \cite{Miller:2003,pz1}, the duality issues of linear programs were studied in \cite{helmes,pz2}, and numerical methods for infinite linear programs were developed in \cite{cho,gai}.
\item For application of the impulse control theory to, e.g., epidemiology, nanoelectronics, Internet congestion control, medicine and inventory, see \cite{AvrachenkovHabachiPZ:2015,helmes,hou,PiunovskiySasha:2018,ppt,yang}.
\end{itemize}

The rest of this article is organized as follows. After the preliminary model description in Section \ref{sec2}, we provide the rigorous problem statement which is most similar to the problems in \cite{chris,pz1,pz2}. Note that there were no functional constraints in \cite{chris,pz2} and no discounting in \cite{pz1,pz2}. The method of investigation, dynamic programming, was also different in \cite{chris}; duality was not studied in \cite{chris,pz1}. All this confirms the novelty of the current work.  We reformulate the impulse control problem as a constrained MDP problem in Section \ref{sec3}. In Section \ref{PZ2025Section01}, we present relevant results for the induced MDP, formulate the relaxed problem, the primal convex program and its dual program, and observe the absence of relaxation gap and the duality gap.  The main results are presented in Section \ref{PZ2025Section02}, consisting of the formulation of the second convex program and its dual program, the duality results, and the general procedure for finding optimal  $(J+1)$-mixed strategies for the original impulse control problems.  This procedure is illustrated with a solved example in Section \ref{sec7}. Section \ref{secap} contains the proofs of the statements.

The following notations are used throughout this paper. Measures are  $[0,\infty]$-valued. We use the notations $\RR_+:=[0,\infty)$ and $\RR_+^{J}:=[0,\infty)^J$. When  ${\bf Y}$ is considered as a subset of ${\bf X}$, the complement of the set ${\bf Y}$ in ${\bf X}$ is denoted by ${\bf Y}^c.$ Let ${\bf V}\subseteq {\bf Y}$ be fixed. The restriction on ${\bf V}$ of a $[-\infty,\infty]$-valued function $f(\cdot)$ defined on ${\bf Y}$ is denoted by $f|_{\bf V}$. This convention applies to measures as functions defined on $\sigma$-algebras, too. If $({\bf Y},{\cal F})$ and $({\bf B},{\cal G})$ are measurable spaces, and  $\mu$ is a measure on ${\bf Y}\times {\bf B}$ endowed with the product $\sigma$-algebra, then the marginal of $\mu$ on ${\bf Y}$ is denoted by $\mu_{\bf Y}$ defined by  $\mu_{\bf Y}(dy)=\mu(dy\times {\bf B}).$  We denote by ${\cal M}^{\sigma<\infty}({\bf Y})$ (or ${\cal M}^{<\infty}({\bf Y})$) the set of $\sigma$-finite (resp., finite) measures on $({\bf Y},{\cal F})$ or say for brevity on ${\bf Y}$. The space of all bounded measurable functions on ${\bf Y}$ is denoted by $\mathbb{B}({\bf Y}).$

\section{Preliminary model description}\label{sec2}
The state space  ${\bf X}$ of the system under investigation is a non-empty Borel subset of a complete separable metric space (i.e., Polish space)  with metric $\rho_X$ and the Borel $\sigma$-algebra ${\cal B}({\bf X})$. Typical examples are the finite-dimensional Euclidean spaces $\RR^n$. In absence of control, the dynamics is described by the (semi-)flow $\phi(\cdot,\cdot):~{\bf X}\times[0,\infty)\to{\bf X}$.   This flow possesses the following properties:
\begin{itemize}
\item[(a)] $\phi(x,0)=x$ for all $x\in{\bf X}$;
\item[(b)] the mapping $\phi(\cdot,\cdot)$ is measurable, where and below measurability is understood as Borel-measurability;
\item[(c)] $\lim_{t\downarrow 0}\phi(\phi(x,s),t)=\phi(x,s)$ for all $x\in{\bf X}$ and $s\in[0,\infty)$;
\item[(d)] $\phi(x,(s+t))=\phi(\phi(x,s),t)$ for all $x\in{\bf X}$ and $(s,t)\in[0,\infty)\times[0,\infty)$.
\end{itemize}
Suppose at the initial time moment $0$ that the state of the system is $x_0\in{\bf X}$. As time goes on, in absence of impulses, the process moves along the flow $\phi$ leading to the  trajectory $\{\phi(x_0,t),~t\ge 0\}$. Property (c) of the flow $\phi$ means that the trajectories of the (uncontrolled) system are continuous from the right, whereas property (d) means that, in absence of impulses, the movement along the flow from $x\in{\bf X}$ to $\phi(x,s)$ on the time interval $(0,s]$ and the further movement from $\phi(x,s)$ to $\phi(\phi(x,s),t)$ on the time interval $(s,s+t]$ can be equivalently represented as the combined movement from $x$ to $\phi(x,s+t)$ on the interval $(0,s+t]$.

\begin{figure}[htbp]
\begin{center}
\includegraphics[width=11cm]{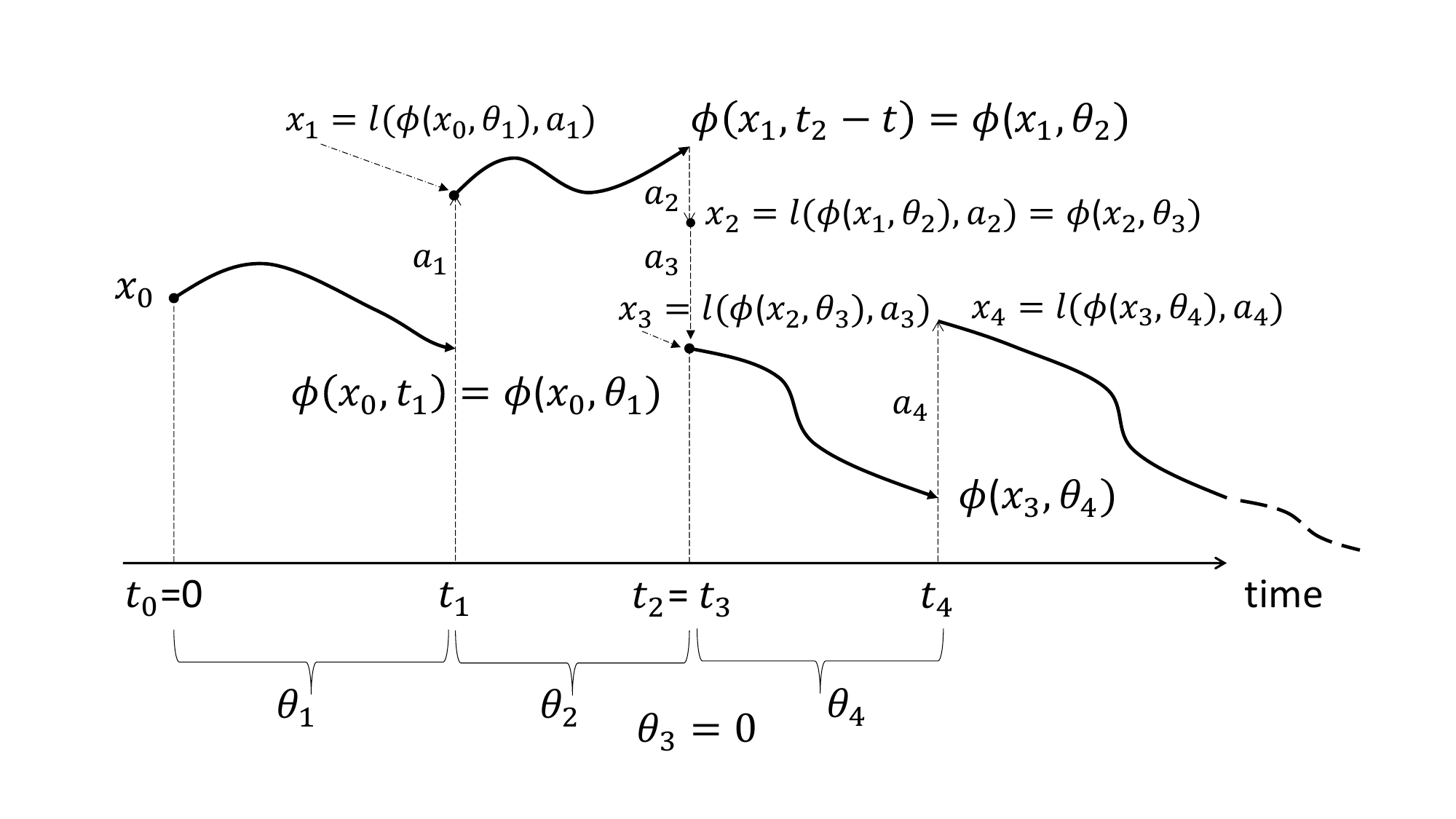}
\caption{Trajectories of the impulsively controlled dynamical system.}\label{figure111}
\end{center}
\end{figure}

Roughly speaking, in the system under (purely) impulsive control, the decision maker decides particular time instants $0\le t_1\le t_2\le\dots\le\infty$, at which he applies impulsive controls (or say impulses) $a_1,a_2,\ldots\in {\bf A}$ with $\bf A$ being a non-empty Borel subset of a complete separable metric space with metric $\rho_A$ and the Borel $\sigma$-algebra ${\cal B}({\bf A})$. Here, if $t_n=\infty$, the impulse applied at $t_n$ is fictitious, and formally speaking, it means no further impulses are applied after $t_{n-1}$.
Each impulse $a\in{\bf A}$ applied at the current state $x\in{\bf X}$ results in the immediate jump of the process to the new state $l(x,a)\in{\bf X}$, where the mapping $l(\cdot,\cdot):~{\bf X}\times{\bf A}\to{\bf X}$ is assumed to be measurable. The case when $l(x,a)=x$ for some $a\in{\bf A}$ is not excluded. After an impulse say $a$ is applied at $t_n<\infty$ when the state is $x$, and before the next time instant $t_{n+1}$, the system dynamics evolves according to the flow $\phi(l(x,a),t)$ with $t$ being the time elapsed from the last impulse time instant $t_n$. If two finite consecutive impulse time instants coincide, say $t_{n}=t_{n+1}<\infty$, then the state resulted in by the impulse $a'$ at $t_{n+1}$ is $l(l(x,a),a')$. This happens when multiple impulses are applied at a single time moment (see $t_3=t_2$ in Fig.\ref{figure111}), and we do not exclude such situations from consideration. It implies that at the time moment $t_n$, multiple values of the state are assigned. While it is possible to describe such multi-valued trajectories of the system under impulsive control, for our purpose it is equivalent and more convenient to consider this system as one in discrete-time, where each step corresponds to a time instant when an impulse is applied. This becomes clearer after we describe the performance measure in the concerned optimal control problem. Further details are given in the next section when this discrete-time system, or more precisely, this Markov decision process (MDP), is formulated rigorously.

Now we describe the optimal impulse control problem.  Assume that in absence of impulses, there is a running cost accumulated at the (gradual) cost rate $C^g(x)\in[0,+\infty]$ associated with any state $x\in{\bf X}$. The function $C^g(\cdot)$ is again assumed to be measurable. Besides, whenever an impulse $a$ is applied at state $x$ at some finite time instant, the lump sum cost $C^I(x,a)\in[0,+\infty]$ is incurred, where the function $C^I(\cdot,\cdot)$ is again measurable. See Fig.\ref{figure111}. On top of that, we assume that the costs (running cost and impulse cost) are discounted with the discount factor $\alpha\in(0,\infty).$
In the current section, to simplify the presentation, we consider only one performance measure coming from the functions $C^g(\cdot)$ and $C^I(\cdot,\cdot)$, see the next paragraph.

 It is convenient to introduce $\theta_i:=t_i-t_{i-1},~i=1,2,\ldots$, accepting $\infty-\infty=\infty,$ $t_0:=0$, so that $t_i=\sum_{k=1}^i\theta_k,~i=0,1,2,\ldots$. In terms of interpretation, when it is finite, $\theta_i$ is the time  duration between two corresponding consecutive time instants when impulses are applied. If, with some $N\in\{1,2,\dots\}$,  there are $N$ impulses $a_1,a_2,\dots,a_N$ being applied at the time moments $t_1,t_2,\ldots,t_N$, respectively, and no further impulses beyond $t_N$, meaning that $\theta_{N+1}=t_{N+1}=t_{N+2}=\dots=\infty$,  then the total (discounted) cost equals
\begin{eqnarray}
&&\sum_{i=1}^N\left\{\int_0^{\theta_i} e^{-\alpha(t_{i-1}+t)} C^g(\phi(x_{i-1},t))~dt+e^{-\alpha t_i}C^I(\phi(x_{i-1},\theta_i),a_i)\right\}\nonumber\\
&&+\II\{N<\infty\}\int_0^\infty e^{-\alpha(t_{N}+t)} C^g(\phi(x_N,t))~dt,\label{eq1}
\end{eqnarray}
where $\alpha>0$ is the fixed discount factor.
In the above formula, $x_i:=l(\phi(x_{i-1},\theta_i),a_i)$, $i=1,2,\ldots,N$. Expression (\ref{eq1}) is the performance measure of a particular (control) strategy applying those specified impulses at the given finitely many time instants. Note that the indicator $\II\{N<\infty\}$ was redundant because $N<\infty$ under the current strategy. Nevertheless, we included it to make the formula valid even when $N=\infty$, provided that we accept that $\II\{N<\infty\}\int_0^\infty e^{-\alpha(t_{N}+t)} C^g(\phi(x_N,t))~dt=0$  when $N=\infty,$ and $e^{-\alpha t_i}C^I(\phi(x_{i-1},\theta_i),a_i)=0$ when $t_i=\infty.$  In general, the number of impulses is not bounded {\it a priori}, neither are all $t_i$ (or $a_i$) necessarily deterministic time instants (resp., impulses). They depend on the control strategy. Accordingly, the expression of the concerned performance measures of a general strategy will be given in the next section after the definitions of a strategy and its corresponding so-called strategic measure become clear. For now, one may keep in mind that the optimal control problem investigated in this paper is about minimizing over strategies the expected total (discounted) cost, subject to further constraints to be introduced.

\section{Reformulation as an MDP}\label{sec3}
We now formulate the concerned optimal impulse control problem as an MDP.
As mentioned in the last section, this MDP naturally arises from inspecting the system state $\{x_{i-1}\}_{i=1}^\infty$ after each impulse is applied, after the initial state $x_0$ is observed. Certainly, the sequences $\{\theta_i\}_{i=1}^\infty$ and $\{t_i\}_{i=1}^\infty$ are in one-to-one correspondence with $t_i=\sum_{j=1}^i \theta_j$. The actions in this MDP are the pairs of $(\theta_i,a_i)$ with $\theta_i$ being the time to wait until the next time instant when the impulse $a_i$ is applied. Then in this MDP, the one-step cost is formed by the (discounted) running cost as well as the impulse cost incurred over two consecutive time instants when impulses are applied in the form of
\begin{eqnarray*}
&&e^{-\alpha t_{i-1}}\int_0^{\theta_i} e^{-\alpha t} C^g(\phi(x_{i-1},t))~dt+e^{-\alpha t_i}C^I(\phi(x_{i-1},\theta_i),a_i)\\
&=&e^{-\alpha t_{i-1}}\left\{\int_0^{\theta_i} e^{-\alpha t} C^g(\phi(x_{i-1},t))~dt+e^{-\alpha \theta_i}C^I(\phi(x_{i-1},\theta_i),a_i)\right\}.
\end{eqnarray*}
Here we accept that, if $\theta=\infty$, then $e^{-\alpha\theta}=0$ and $0 \times \infty=0.$
However, if one pursues this way, then due to the presence of discounting (more precisely, $e^{-\alpha t_{i-1}}$) in the above expression, it would enforce the state of the MDP at step $i-1$ to be in form of $(t_{i-1},x_{i-1})$. For brevity, we would prefer to avoid the inclusion of the extra coordinate $t_{i-1}$ into the state of the MDP. This is done by observing that $e^{-\alpha t_{i-1}}=\prod_{j=1}^{i-1} e^{-\alpha \theta_j}$ and $e^{-\alpha \theta_j}\in [0,1]$ can be understood as the conditional probabilities of certain events described as follows. After an impulse $a_{j}$ is applied at $t_j=t_{j-1}+\theta_j$, the process is artificially terminated (or say ``killed'') with the conditional probability $1-e^{-\alpha \theta_j}$ given that the process has not been killed before. (If $t_{j-1}<\infty$ and $\theta_j=\infty$, it means that no further impulses will be applied beyond $t_{j-1}$, and in this case, the process is killed with the conditional probability $1-e^{-\alpha \infty}=1.$) This amounts to sending the state to $\Delta$, a costless cemetery point, not in ${\bf X}.$ We join $\Delta$ to ${\bf X}$, and form $X_{\Delta}:={\bf X}\cup\{\Delta\}$, where $\Delta$ is an isolated point. $\Delta$ is absorbing in the sense that once $\Delta$ is hit, the system state remains there forever. With this interpretation, one can now take the one-step cost function in the MDP given by \begin{eqnarray*}
\II\{x_{i-1}\in {\bf X}\}\left\{\int_0^{\theta_i} e^{-\alpha t} C^g(\phi(x_{i-1},t))~dt+e^{-\alpha\theta_i}C^I(\phi(x_{i-1},\theta_i),a_i)\right\}:
\end{eqnarray*}
the term $e^{-\alpha t_{i-1}}$ has been replaced with $\II\{x_{i-1}\in {\bf X}\}$, after we extend the state space to ${\bf X}_\Delta$ and modify the transition law  in line with what was described. The resulting MDP is what we deal with throughout the rest of this paper.

The rigorous formulation of the MDP model is as follows. It is given
by the following elements $\{{\bf X}_\Delta, {\bf B}, Q, \{\bar{C}_j\}_{j=0}^J,\{d_j\}_{j=1}^J\}$, where
\begin{itemize}
\item ${\bf X}_\Delta:={\bf X}\cup\{\Delta\}$ is the state space;
\item ${\bf B}:=[0,\infty]\times{\bf A}$ is the action space equipped with the product $\sigma$-algebra, where $[0,\infty]$ is the one-point compactification of $[0,\infty)$;
\item $Q(dy|x,b=(\theta,a)):=\left\{\begin{array}{ll}
e^{-\alpha\theta}\delta_{l(\phi(x,\theta),a)}(dy)+(1-e^{-\alpha\theta})\delta_\Delta (dy), & \mbox{if } x\ne\Delta,\theta\ne\infty;\\
\delta_\Delta(dy) & \mbox{otherwise} \end{array}\right.$ is the transition law;
\item For each $j\in\{0,\dots,J\}$ with $J$ being a nonnegative integer,
    \begin{eqnarray*}
    \bar {C}_j(x,(\theta,a)):=\left\{\begin{array}{ll}
\int_0^\theta e^{-\alpha t}C^g_j(\phi(x,t))dt
+e^{-\alpha\theta}C^I_j(\phi(x,\theta),a), & \mbox{ if } x\in{\bf X};\\
0, & \mbox{ if } x=\Delta,\end{array}\right.
\end{eqnarray*}
where $C^g_j(\cdot)$ and $C^I_j(\cdot,\cdot)$ are $[0,\infty]$-valued measurable functions on ${\bf X}$ and ${\bf X}\times {\bf A}$, respectively, and we accept that $e^{-\alpha\theta}C^I_j(\phi(x,\theta),a)=0$ when $\theta=\infty$;
\item For each $1\le j\le J$, $d_j$ is a nonnegative integer.
\end{itemize}
Note that in the above MDP model, instead of a single cost function $\bar{C}(\cdot,\cdot)$, we have $J+1$ cost functions. This is because we shall consider a constrained optimal impulse control problem in what follows. In this regard, $J$ is the number of constraints, and $\{d_j\}_{j=1}^J$ are the corresponding constraint constants.

Consider the MDP model $\{{\bf X}_\Delta, {\bf B}, Q, \{\bar{C}_j\}_{j=0}^J,\{d_j\}_{j=1}^J\}$.  As usual, we introduce canonical sample space as the countable product $\Omega:=({\bf X}_\Delta\times{\bf B})^\infty$ endowed with the product topology, and the corresponding product $\sigma$-algebra $\cal F$. It coincides with the Borel $\sigma$-algebra on $\Omega.$ An element
$\omega
=(x_0,b_1=(\theta_1,a_1),x_1,b_2=(\theta_2,a_2),x_2,\ldots)\in \Omega
$  is
a trajectory in of the MDP. Let, for each $i\in\{0,1,\dots\}, $
\begin{eqnarray*}
X_i(\omega):=x_i,B_{i+1}(\omega):=b_{i+1},\Theta_{i+1}(\omega):=\theta_{i+1}, A_{i+1}(\omega):=a_{i+1}
\end{eqnarray*} 
be the corresponding projection mappings of $\omega\in\Omega$, and let $
T_i(\omega):=\sum_{j=1}^i \Theta_j(\omega).
$
The process $\{(X_{i-1},B_i)\}_{i\ge 1}$ is often called an MDP, too, but we use the term for both the bivariate process and the underlying model when there is no danger of confusion.
In general, capital letters (not in special fonts) are used to denote random elements defined on $\Omega$, and the argument $\omega$ is typically omitted in expressions inside a probability or an expectation.

Finite sequences of the form
$$h_i=(x_0,b_1=(\theta_1,a_1),x_1,b_2=(\theta_2,a_2),\ldots,x_i)$$
will be called $i$-histories, $i=0,1,\dots$. The space of all such $i$-histories will be denoted as ${\bf H}_i$. We endow it with the Borel $\sigma$-algebra ${\cal F}_i:={\cal B}({\bf H}_i)$.

A (control) strategy $\pi=\{\pi_i\}_{i=1}^\infty$ is a sequence of stochastic kernels $\pi_i$ on $\bf B$ given $\textbf{H}_{i-1}$.
A strategy $\pi$ is called Markov and denoted by $\pi^M$ if for each $i\in\{1,2,\dots\}$, there is some stochastic kernel $\pi_i^M$ on ${\bf B}$ given ${\bf X}_\Delta$ such that $\pi_i(db|h_{i-1})=\pi_i^M(db|x_{i-1})$ for all $h_{i-1}=(x_0,b_1=(\theta_1,a_1),x_1,b_2=(\theta_2,a_2),\ldots,x_{i-1})\in {\bf H}_{i-1}.$
 A Markov strategy $\pi^M$ is called stationary and denoted as ${\pi}^s$, if for all $i\in\{1,2,\dots\}$, there is a stochastic kernel ${\pi}^s$ on $\bf B$ given $\textbf{X}_\Delta$  such that $\pi_i^M(db|x)={\pi}^s(db|x)$ for all $x\in {\bf X}_{\Delta}.$
 A strategy $\pi$ (or Markov strategy $\pi^M$) is called deterministic (resp., deterministic Markov),  if for all $i\in\{1,2,\dots\}$, there is some measurable mapping $f_i(\cdot):\textbf{X}_\Delta\to{\bf B}$ (resp., $f_i^M(\cdot): {\bf X}_{\Delta}\to {\bf B}$) such that  $\pi_i(db|h_{i-1})=\delta_{f_i(x_{i-1})}(db)$ (resp., $\pi_i^M(db|x_{i-1})=\delta_{f_i^M(x_{i-1})}(db)$). Finally, a stationary strategy $\pi^s$ is called deterministic stationary and denoted as $f$,  if $\pi^s(db|x)=\delta_{f(x)}(db)$ for all $x\in {\bf X}_{\Delta}$, where $f(\cdot):\textbf{X}_\Delta\to{\bf B}$ is a measurable mapping.

In terms of interpretation, under a strategy $\pi$, $\pi_i(db|h_{i-1})$ is the (regular) conditional distribution of $B_i=(\Theta_{i},A_i)$ given the $(i-1)$-history $H_{i-1}=h_{i-1}$. This is in line with the following construction of the underlying probability.

For a given initial state $x_0\in \textbf{X}_\Delta$ and a strategy $\pi$, by the Ionescu-Tulcea theorem, see e.g., \cite[Prop.7.28]{Bertsekas:1978}, there is a unique probability measure $\PP^\pi_{x_0}$ on $(\Omega,{\cal F})$ satisfying the following conditions:
\begin{eqnarray*}
\PP^\pi_{x_0}(X_0\in\Gamma_X)=\delta_{x_0}(\Gamma_X)~\forall~ \Gamma_X\in{\cal B}({\bf X}_\Delta),
\end{eqnarray*}
and for all $i=1,2,\dots$, $\Gamma_B\in{\cal B}({\bf B})$, $\Gamma_X\in{\cal B}({\bf X}_\Delta)$, $\PP^\pi_{\nu}$-almost surely
\begin{eqnarray}\label{e33}
&&\PP^\pi_{x_0}(B_i\in\Gamma_B|H_{i-1})=\pi_i(\Gamma_B|H_{i-1});\nonumber\\
&&\PP^\pi_{x_0}(X_i\in\Gamma_X|H_{i-1},B_i=(\Theta_i,A_i))= Q(\Gamma_X| X_{i-1},B_i)  \nonumber\\
&=&\left\{\begin{array}{ll} e^{-\alpha\Theta_i}
\delta_{l(\phi(X_{i-1},\Theta_i),A_i)}(\Gamma_X)\nonumber\\
+(1-e^{-\alpha\Theta_i})\delta_\Delta(\Gamma_X), & \mbox{ if } X_{i-1}\ne\Delta,~\Theta_i\ne\infty; \nonumber\\
\delta_\Delta(\Gamma_X) & \mbox{ otherwise. } \end{array} \right. \\
\end{eqnarray}
The probability $\PP^\pi_{x_0}$ is called the strategic measure of $\pi$ (given the initial state $x_0$). The mathematical expectation with respect to $\PP^\pi_{x_0}$ is denoted as $\EE^\pi_{x_0}$, respectively. Note that $\PP^\pi_{x_0}(X_i\in{\bf X})=\EE^\pi_{x_0}[e^{-\alpha T_i}]$, $i=0,1,\ldots$.

Now we define the performance measures of each strategy.
It will be convenient to introduce the following notation of the performance measure of a generic cost function $\bar{C}(\cdot,\cdot)$ on ${\bf X}_{\Delta}\times {\bf B}$ for the MDP in the form of
   \begin{eqnarray*}
    \bar {C}(x,(\theta,a)):=\left\{\begin{array}{ll}
\int_0^\theta e^{-\alpha t}C^g(\phi(x,t))dt
+e^{-\alpha\theta}C^I(\phi(x,\theta),a), & \mbox{ if } x\in{\bf X};\\
0, & \mbox{ if } x=\Delta,\end{array}\right.
\end{eqnarray*}
with $C^I(\cdot,\cdot)$ and $C^g(\cdot)$ being  $[0,\infty]$-valued  measurable functions on ${\bf X}\times {\bf A}$ and on ${\bf X}$, respectively. That $\bar{C}(\cdot,\cdot)$ is called generic in the sense that if $C^I(\cdot,\cdot)$ and $C^g(\cdot)$ in its definition are replaced by $C^I_j(\cdot,\cdot)$ and $C^g_j(\cdot)$, then we obtain $\bar{C}_j(\cdot,\cdot)$. Now we define for each $x\in {\bf X}_\Delta$ and strategy $\pi$  the performance measure
\begin{eqnarray*}
{\cal V}(x,\pi,\bar{C}) :=\EE^\pi_x\left[\sum_{i=1}^\infty\bar C(X_{i-1},(\Theta_i,A_i))\right].
\end{eqnarray*}

Below, we take $x_0\in\textbf{X}$ as a fixed initial state, since the case of $x_0=\Delta$ is trivial. Then the constrained optimal control problem under study is the following one:
\begin{eqnarray}\label{PZZeqn02}
\mbox{Minimize with respect to } \pi &&{\cal V}(x_0, \pi,\bar{C}_0)  \\
\mbox{subject to }&& {\cal V}(x_0,\pi,\bar{C}_j)\le d_j,~j=1,2,\dots,J.\nonumber
\end{eqnarray}
Note that we do not exclude the possibility of $\sum_{i=1}^{\infty}\Theta_i<\infty$, but we will only consider the total cost accumulated over $[0,\sum_{i=1}^\infty \Theta_i).$ This is consistent with the definition of ${\cal V}(x_0, \pi, \bar{C}_j)$.

\begin{definition}\label{d2}
A strategy $\pi$ is called  feasible in problem (\ref{PZZeqn02}) if it satisfies all the constraint inequalities therein. A strategy $\pi^\ast$ is called optimal in problem (\ref{PZZeqn02}) if it is feasible and ${\cal V}(x_0,\pi^*,\bar{C}_0)\le  {\cal V}(x_0,\pi,\bar{C}_0)$ for all feasible strategies $\pi$ in problem (\ref{PZZeqn02}).  A feasible strategy $\pi$ such that ${\cal V}(x_0,\pi,\bar{C}_0)<\infty$ is called feasible with a finite value.
\end{definition}

The infimum of  ${\cal V}(x_0,\pi,\bar C_0)$ over all feasible strategies is called the  optimal value of problem (\ref{PZZeqn02}), and is denoted by $val(\mbox{problem (\ref{PZZeqn02})})$. The notation $val(\cdot)$ will also be used for other concerned optimization problems; its meaning lies unambiguously in its argument.
Note that, if there are no feasible strategies with a finite value, then 
 either problem (\ref{PZZeqn02}) is not consistent or any feasible strategy is optimal.

\section{Facts and preliminaries}\label{PZ2025Section01}

\subsection{Semicontinuous induced MDP}

The MDP (model) $\{{\bf X}_\Delta, {\bf B}, Q, \{\bar{C}_j\}_{j=0}^J,\{d_j\}_{j=1}^J\}$ introduced earlier in the previous section is induced by the impulse control model, and thus will be referred to as the induced MDP. If the state space ${\bf X}_\Delta$ is replaced with an arbitrary nonempty Borel space say ${\bf Y}$,   and ${\bf B}$, $Q$ and $\{\bar{C}_j\}_{j=0}^J$ are arbitrary nonempty Borel action space, stochastic kernel and $[0,\infty]$-valued cost functions (not necessarily induced in the above specific form by $C^g_j(\cdot)$ and $C^I_j(\cdot,\cdot)$), then we call the MDP $\{{\bf Y},{\bf B}, Q, \{\bar{C}_j\}_{j=0}^J,\{d_j\}_{j=1}^J\}$ generic. The optimal control problem (\ref{PZZeqn02}), Definition \ref{d2} and Definition \ref{15JuneDef02} below can be understood for a generic MDP, too, with obvious adjustments; in particular, $x_0$ will be a fixed state in ${\bf Y}$.

\begin{definition}\label{PZ2025Def01}
A generic MDP $\{{\bf Y},{\bf B}, Q, \{\bar{C}_j\}_{j=0}^J,\{d_j\}_{j=1}^J\}$ is called semicontinuous if the following hold:
\begin{itemize}
\item[(a)] the action space ${\bf B}$ is compact;
\item[(b)] for each bounded continuous function ${u}(\cdot)$ on ${\bf Y}$, $\int_{\bf Y}u(y)Q(dy|x,b)$ is continuous in $(x,b)\in{\bf Y}\times{\bf B}$;
    \item[(c)] for each $j\in\{0,\dots,J\},$  $\bar{C}_j(\cdot,\cdot)$ is lower semicontinuous.
\end{itemize}
\end{definition}

The following condition implies that the induced MDP is semicontinuous.

\begin{condition}\label{ConstrainedPPZcondition01}
\begin{itemize}
\item[(a)]  $\textbf{A}$ is compact.
\item[(b)] $l(\cdot,\cdot)$ is continuous on ${\bf X}\times {\bf A}$.
\item[(c)] $\phi(\cdot,\cdot)$ is continuous on ${\bf X}\times [0,\infty)$.
\item[(d)] For each $j\in\{0,\dots,J\}$,  $C_j^I(\cdot,\cdot)$ is lower semicontinuous on ${\bf X}\times {\bf A}$, and  $C_j^g(\cdot)$ is lower semicontinuous on ${\bf X}$.
\end{itemize}
\end{condition}

\begin{proposition}\label{pr1}
Under Condition \ref{ConstrainedPPZcondition01}, the induced MDP is semicontinous.
\end{proposition}
The proofs of this proposition and most of the statements are presented in Section \ref{secap}.

\begin{lemma}\label{l1}
If the induced MDP is semicontinuous and admits a feasible strategy with a finite value, then there exists a stationary strategy, which is with a finite value and optimal in problem  (\ref{PZZeqn02}).
\end{lemma}
\par\noindent\textit{Proof}. The assertion follows from the fact established in Theorem 4.1 of \cite{Dufour:2012} for a generic MDP. $\hfill\Box$

Let
\begin{eqnarray}\label{e5p}
{\bf V}:=\left\{x\in\textbf{X}_\Delta: \inf_{\pi}\EE_x^\pi\left[ \sum_{i=1}^\infty \sum_{j=0}^J\bar{C}_j(X_{i-1},B_{i}) \right]=\inf_\pi{\cal V}\left(x,\pi,\sum_{j=0}^J \bar{C}_j\right)>0\right\},
\end{eqnarray}
Note that $\Delta\in {\bf V}^c:=\textbf{X}_\Delta\setminus {\bf V}$ because $\Delta$ is a costless cemetery in the induced MDP.

Suppose the induced MDP is semicontinuous. Applying Corollary 9.17.2 of \cite{Bertsekas:1978} to the induced MDP under study with the cost function $\sum_{j=0}^J \bar C_j(\cdot,\cdot)$, we obtain that the  Bellman function 
\begin{equation}\label{e8}
{\cal V}^*(x,\sum_{j=0}^J\bar C_j):=\inf_\pi \EE^\pi_x\left[\sum_{i=1}^\infty \sum_{j=0}^J \bar C_j(X_{i-1},B_i)\right],~x\in{\bf X}_\Delta
\end{equation}
is a non-negative  lower semicontinuous function, so that, in particular, the set ${\bf V}^c$ is a closed subset of ${\bf X}_\Delta$ and $\bf V$ is an open subset of $\bf X$. Moreover, there is a deterministic stationary strategy $f^*$ providing the infimum in (\ref{e8}) and in the Bellman equation satisfied  by ${\cal V}^*(\cdot,\sum_{j=0}^J\bar C_j)$:
$${\cal V}^*(x,\sum_{j=0}^J\bar C_j)=\inf_{b\in{\bf B}}\left\{\sum_{j=0}^J \bar C_j(x,b)+\int_{{\bf X}_\Delta}{\cal V}^*(y,\sum_{j=0}^J\bar C_j)Q(dy|x,b)\right\},~x\in{\bf X}_\Delta,$$
and
\begin{equation}\label{e5}
\bar C_j(x,f^*(x))=Q({\bf V}|x,f^*(x))=0~ \forall ~x\in{\bf V}^c, j\in \{0,\dots, J\}.
\end{equation}
For the last assertion, it is sufficient to notice that for all $x\in{\bf V}^c$
$$0={\cal V}^*(x,\sum_{j=0}^J\bar C_j)=\sum_{j=0}^J\bar C_j(x,f^*(x))+\int_{\bf V} {\cal V}^*(y,\sum_{j=0}^J\bar C_j)Q(dy|x,f^*(x))$$
and ${\cal V}^*(y,\sum_{j=0}^J\bar C_j)>0$ for $y\in{\bf V}$. We fix this deterministic stationary strategy $f^\ast.$

Now it is clear from (\ref{e5}) that, if $x_0\in{\bf V}^c$, then $f^*$ is an optimal strategy in problem  (\ref{PZZeqn02}) because
$${\cal V}(x,f^*,\bar{C}_j)=0~\forall~ x\in{\bf V}^c,~~j=0,1,\ldots, J.$$ Consequently, in what follows, we suppose that $x_0\in {\bf V}$, and thus ${\bf V}\ne\emptyset.$
It is also clear that, for each $i\in\{1,2,\dots\}$, it is optimal to apply actions $f^*(X_{i-1})$ as soon as $X_{i-1}\in{\bf V}^c$ because in this way the expected subsequent costs associated with $\bar C_j(\cdot,\cdot)$, $j=0,1,\ldots$, would be zero.

\begin{definition}\label{d3}
Suppose that the induced MDP is semicontinuous and admits a  feasible strategy with a finite value. We denote by
$\Pi$ the set of feasible strategies with a finite value that satisfy
\begin{equation}\label{e0}
\pi_i(db|h_{i-1}):=\delta_{f^*(x_{i-1})}(db),~\forall~ i=1,2,\ldots,
\end{equation}
whenever $x_{i-1}$, the last component of $h_{i-1}$, is in ${\bf V}^c$.
\end{definition}
Now, given that the induced MDP is semicontinuous and admits a feasible strategy with a finite value, it is sufficient to restrict to the strategies from $\Pi$ and  concentrate on the subset $\bf V$. After some obvious modifications if necessary, one can regard, without loss of generality, that the optimal strategy $\pi$ in problem (\ref{PZZeqn02}), existing by Lemma \ref{l1}, is from the class $\Pi$.

\subsection{Relaxed problem and absence of relaxation gap}

\begin{definition}\label{15JuneDef02}
For each strategy $\pi$ in the induced MDP, its  occupation measure $\mu^\pi$  is defined by
$$\mu^\pi(\Gamma):=\EE_{x_0}^\pi\left[\sum_{i=0}^\infty \II\{(X_i,B_{i+1})\in\Gamma\}\right]~\forall~\Gamma\in{\cal B}(\textbf{X}_\Delta\times \textbf{B}).$$
\end{definition}

It holds for each $j\in\{0,\dots,J\}$ that 
$${\cal V}(x_0,\pi,\bar{C}_j)=\int_{{\bf X}_\Delta\times{\bf B}} \bar C_j(x,b)\mu^\pi(dx\times db), ~j=0,1,\ldots, J.$$
Hence, problem (\ref{PZZeqn02}) is equivalent to minimizing $\int_{{\bf X}_\Delta\times{\bf B}} \bar C_j(x,b)\mu^\pi(dx\times db)$ over the space of occupation measures of feasible strategies in the induced MDP. As is well known (see Lemma 9.4.3 of \cite{HernandezLerma:1999}), each occupation measure $\mu$ satisfies equation
\begin{eqnarray}\label{enum1}
\mu(\Gamma\times \textbf{B})=\delta_{x_0}(\Gamma)+\int_{\textbf{X}_\Delta\times \textbf{B}}Q(\Gamma|y,b)\mu(dy\times db),~~~\Gamma\in{\cal B}(\textbf{X}_\Delta).
\end{eqnarray}

\begin{lemma}\label{l2}
Suppose the induced MDP is semicontinuous. Then for each measure $\mu$ on ${\bf X}_\Delta\times{\bf B}$, satisfying equation (\ref{enum1}) and such that
$  \int_{{\bf X}_\Delta\times{\bf B}}\sum_{j=0}^J \bar C_j(x,b)\mu(dx\times db)<\infty,$ it holds that 
$\mu_{{\bf X}_\Delta}|_{\bf V}\in {\cal M}^{\sigma<\infty}({\bf V}),$
where ${\cal M}^{\sigma<\infty}({\bf V})$ is the set of $\sigma$-finite measures on $({\bf V},{\cal B}(\bf V)).$ In particular, for each feasible strategy $\pi$ with a finite value, it holds that $\mu^\pi_{{\bf X}_\Delta}|_{{\bf V}}\in {\cal M}^{\sigma<\infty}({\bf V}).$
\end{lemma}
For the proof of the first assertion, see Theorem 3.2 of \cite{Dufour:2012}. The last assertion follows from the first one, because for each feasible strategy $\pi$ with a finite value, its occupation measure $\mu^\pi$ satisfies (\ref{enum1}) and $\int_{{\bf X}_\Delta\times{\bf B}}\sum_{j=0}^J \bar C_j(x,b)\mu^\pi(dx\times db)=\sum_{j=0}^J{\cal V}(x_0,\pi,\bar{C_j})<\infty.$

Suppose the induced MDP is semicontinuous and admits a feasible strategy with a finite value. According to the explanations at the end of the previous subsection (see (\ref{e5})), for  problem (\ref{PZZeqn02}) it is sufficient to restrict to the set of strategies $\Pi\ne \emptyset$. That, in turn, is equivalent to  minimizing $\int_{{\bf X}_\Delta\times{\bf B}} \bar C_j(x,b)\mu^\pi(dx\times db)$ over the space of occupation measures of feasible strategies in the induced MDP. We shall introduce a relaxation of the latter problem. To this end, let us now examine the occupation measures of strategies in $\Pi$ and make two observations.

\begin{itemize}
\item[($\alpha$)] If $\pi\in\Pi$, then for each measurable subset $\Gamma_1\subset {\bf V}^c$,
$$\mu^\pi(\Gamma_1\times\Gamma_2)=\EE_{x_0}^\pi\left[\sum_{i=0}^\infty \II\{(X_i \in\Gamma_1\}\II\{\Gamma_2\ni f^*(X_i)\}\right]~\forall~\Gamma_2\in{\cal B}(\textbf{B}).$$
Therefore,
\begin{equation}\label{e6}
\int_{{\bf V}^c\times{\bf B}} \bar C_j(x,b)\mu^\pi(dx\times db)=\EE^\pi_{x_0}\left[\sum_{i=0}^\infty \II\{X_i\in{\bf V}^c\} \bar C_j(X_i,f^*(X_i))\right]=0,~
 j=0,1,\ldots, J
\end{equation}
by (\ref{e5}), and
$${\cal V}(x_0,\pi,\bar{C}_j)=\int_{{\bf V}\times{\bf B}} \bar C_j(x,b)\mu^\pi(dx\times db), ~j=0,1,\ldots, J.$$
\item[($\beta$)] For $\pi\in\Pi$,  the measure $\mu=\mu^\pi|_{{\bf V}\times {\bf B}}$ on ${\bf V}\times{\bf B}$ satisfies the equation
$$\mu(\Gamma\times \textbf{B})=\delta_{x_0}(\Gamma)+\int_{\textbf{V}\times \textbf{B}}Q(\Gamma|y,b)\mu(dy\times db),~~\Gamma\in{\cal B}({\bf V}),$$
because for each measurable subset  $\Gamma\subset{\bf V}$,  again by (\ref{e5}),
$$\int_{\textbf{V}^c\times \textbf{B}}Q(\Gamma|y,b)\mu^\pi(dy\times db)=\EE^\pi_{x_0}\left[\sum_{i=0}^\infty \II\{X_i\in{\bf V}^c\} Q(\Gamma|X_i,f^*(X_i))\right]=0.$$
\end{itemize}

We are ready to formulate the following relaxation of problem (\ref{PZZeqn02}) for the induced MDP, assumed to be semicontinuous and to admit a feasible strategy with a finite value: 
\begin{eqnarray}\label{SashaLp02}
\mbox{Minimize}&:&\int_{{\bf V} \times\textbf{B}}\bar{C}_0(x,b)\mu(dx\times db)\\
&&\mbox{ over measures $\mu$ on ${\bf V}\times {\bf B}$ with $\mu_{\bf V}\in {\cal M}^{\sigma<\infty}(\bf V)$} \nonumber  \\
 \mbox{subject to}&:& \mu(\Gamma\times \textbf{B})=\delta_{x_0}(\Gamma)+\int_{{\bf V}\times \textbf{B}}Q(\Gamma|y,b)\mu(dy\times db),~\Gamma\in{\cal B}({\bf V})\label{e9}\\
\mbox{and}&&\int_{{\bf V} \times\textbf{B}}\bar{C}_j(x,b)\mu(dx\times db)\le d_j,~j=1,2,\dots,J. \label{e10}
\end{eqnarray}
The fact that the above problem is a relaxation of the MDP problem (\ref{PZZeqn02}) is because of Items ($\alpha,\beta$) in the above,  Lemma \ref{l2} and the paragraph below it. In particular, 
\begin{eqnarray}\label{PZ2025Eqn01}
val(\mbox{relaxed problem (\ref{SashaLp02})-(\ref{e10})})\le val(\mbox{problem (\ref{PZZeqn02})})<\infty,
\end{eqnarray}
given that  the induced MDP is semicontinuous and admits a feasible strategy with a finite value.
We remark that in general, a measure $\mu$ satisfying equation (\ref{enum1}) or (\ref{e9}) may not be an occupation measure for any strategy $\pi$. As an example, see Subsection  2.3.20 of \cite{ex}. Thus, in general, one cannot conclude that the relaxed problem and the original problem are equivalent \textit {a prioi}. Theorem \ref{t1} below shows that there is no relaxation gap if the induced MDP is semicontinuous and has a feasible strategy with a finite value. It will follow from the next few observations.

\begin{definition}\label{d4}
Suppose the induced MDP is semicontinuous. Consider a  measure $\mu$ on ${\bf V}\times{\bf B}$  with $\mu|_{\bf V}\in {\cal M}^{\sigma<\infty}(\bf V)$. We say that a stationary strategy $\pi^\mu$ is induced by $\mu$,  if it has the following form
$$\pi^\mu(db|x)=\II\{x\in{\bf V}\}\varphi_\mu(db|x)+\II\{x\in{\bf V}^c\}\delta_{f^*(x)}(db),$$
where $\varphi_\mu$ is a stochastic kernel on $\bf B$ given $\bf V$ satisfying
$$\varphi_\mu(db|x)\mu(dx\times{\bf B})=\mu(dx\times db)\mbox{ on } {\cal B}({\bf V}\times{\bf B}).$$
\end{definition}
The stochastic kernel $\varphi_\mu$ exists by Corollary 7.27.1 of \cite{Bertsekas:1978} which should be applied to the sets $\{\Gamma_n\times{\bf B}\}_{n=1}^\infty$ that decompose ${\bf V}\times {\bf B}$ with $\mu(\Gamma_n\times{\bf B})<\infty$. There may be multiple induced strategies by $\mu$, different on  $\mu(dx\times{\bf B})$-null subsets of $\bf V$.

\begin{lemma}\label{l3}
Suppose the induced MDP is semicontinuous. Let a measure $\mu$ be feasible in the   problem  (\ref{SashaLp02})-(\ref{e9}) such that $\int_{{\bf V}\times{\bf B}}\sum_{j=0}^J \bar C_j(x,b)\mu(dx\times db)<\infty$. Then for each $j=0,1,\ldots, J$
$$\EE^{\pi^\mu}_{x_0}\left[\sum_{i=1}^\infty \bar C_j(X_{i-1},B_i)\right]=\int_{{\bf V}\times{\bf B}} \bar C_j(x,b)\mu^{\pi^\mu}(dx\times db)\le\int_{{\bf V}\times{\bf B}} \bar C_j(x,b)\mu(dx\times db).$$
\end{lemma}
For the proof, see Corollary 3.1 of \cite{Dufour:2012}. In that corollary, integration was over ${\bf X}_\Delta\times{\bf B}$, but, similarly to (\ref{e6}), $\displaystyle \int_{{\bf V}^c\times{\bf B}} \bar C_j(x,b)\mu^{\pi^\mu}(dx\times db)=0$. 

Note that under the conditions imposed  in Lemma \ref{l3}, if $\mu$ therein is feasible in the relaxed problem  (\ref{SashaLp02})-(\ref{e10}), then $\Pi$ is non-empty, because the strategy $\pi^\mu$ therein belongs to  $\Pi$.

\begin{theorem}\label{t1}
Suppose the induced MDP is semicontinuous and admits a feasible strategy $\pi$ with a finite value. Then the following assertions hold.
\begin{itemize}
\item[(a)] $val(\mbox{problem (\ref{PZZeqn02})})=val(\mbox{relaxed problem (\ref{SashaLp02})-(\ref{e10})})<\infty$.
\item[(b)] If a strategy $\pi^*\in\Pi$ is optimal in problem (\ref{PZZeqn02}), then the measure $\mu^{\pi^*}|_{{\bf V}\times{\bf B}}$ is optimal in the relaxed problem (\ref{SashaLp02})-(\ref{e10}). Note that an optimal strategy $\pi^*$ exists by Lemma \ref{l1} and explanations at the end of the previous subsection.
\item[(c)] If a measure $\mu^*$ is optimal in the relaxed problem (\ref{SashaLp02})-(\ref{e10}), then the induced strategy $\pi^{\mu^*}\in\Pi$ is optimal in problem (\ref{PZZeqn02}).
\end{itemize}
\end{theorem}

\par\noindent\textit{Proof of Theorem \ref{t1}.} (a) In view of (\ref{PZ2025Eqn01}), we only need to show that $val(\mbox{problem (\ref{PZZeqn02})})$\linebreak$\le val(\mbox{relaxed problem (\ref{SashaLp02})-(\ref{e10})})$ as follows.
Consider an arbitrary measure $\mu$ that is feasible in the relaxed problem  (\ref{SashaLp02})-(\ref{e10}) and satisfies $\int_{{\bf V}\times {\bf B}}\bar{C}_0(x,b)\mu(dx\times db)<\infty$. For the relaxed problem  (\ref{SashaLp02})-(\ref{e10}), it suffices to restrict to the set of such measures, which is non-empty because of (\ref{PZ2025Eqn01}). Then for  the induced (by $\mu$) strategy $\pi^\mu$, 
$${\cal V}(x_0,\pi^\mu,\bar{C}_j)\le \int_{{\bf V}\times{\bf B}} \bar C_j(x,b)\mu(dx\times db),~~j=0,1,\ldots,J$$
by Lemma \ref{l3}. Therefore, the induced strategy $\pi^\mu$ is feasible in problem (\ref{PZZeqn02}) and
$${\cal V}(x_0,\pi^\mu,\bar{C_0})\le \int_{{\bf V}\times{\bf B}} \bar C_0(x,b)\mu(dx\times db).$$ The arbitrariness of $\mu$ and this inequality  lead to  the required inequality.

(b) Let $\pi^*\in\Pi$ be an optimal strategy in problem (\ref{PZZeqn02}). Then the restriction of its occupation measure $\mu^{\pi^*}|_{{\bf V}\times {\bf B}}$ is feasible in the relaxed problem (\ref{SashaLp02})-(\ref{e10}) according to  Lemma \ref{l2} and the observations made below it, see Items ($\alpha,\beta$) therein. Moreover,
\begin{eqnarray*}
{\cal V}(x_0,\pi^*,\bar{C}_0)&=&\int_{{\bf V}\times{\bf B}} \bar C_0(x,b)\mu^{\pi^*}(dx\times db)=\int_{{\bf V}\times {\bf B}}\bar{C}_0(x,b)\mu^{\pi^*}|_{{\bf V}\times {\bf B}}(dx\times db)\\
&=& val(\mbox{problem (\ref{PZZeqn02})})= val(\mbox{relaxed problem (\ref{SashaLp02})-(\ref{e10})}),
\end{eqnarray*}
where the last equality holds by (a). 
Thus, the measure $\mu^{\pi^*}|_{{\bf V}\times {\bf B}}$ is optimal in the relaxed problem (\ref{SashaLp02})-(\ref{e10}).

(c) Let the measure $\mu^*$ be optimal in the relaxed problem (\ref{SashaLp02})-(\ref{e10}). It holds that $\int_{{\bf V}\times{\bf B}} \bar C_0(x,b)$\linebreak$\times\mu^*(dx\times db)<\infty$, because  $val(\mbox{relaxed problem (\ref{SashaLp02})-(\ref{e10})})<\infty$  as was shown in (a). Now, Lemma \ref{l3} is applicable to $\mu^\ast$, and we see that
the induced (by $\mu^\ast$) strategy $\pi^{\mu^*}$  is feasible in problem (\ref{PZZeqn02}), and
\begin{eqnarray*}
{\cal V}(x_0,\pi^{\mu^*},\bar{C}_0)&\le &\int_{{\bf V}\times{\bf B}} \bar C_0(x,b)\mu^*(dx\times db)\\
&=& val(\mbox{relaxed problem (\ref{SashaLp02})-(\ref{e10})})= val(\mbox{problem (\ref{PZZeqn02})}).
\end{eqnarray*}
Thus the strategy $\pi^{\mu^*}$ is optimal in problem (\ref{PZZeqn02}). That $\pi^{\mu^*}$ belongs to $\Pi$ follows from the definitions of $\Pi$ and $\pi^{\mu^\ast}$. \hfill $\Box$

\subsection{Primal convex program and duality results}\label{PZ2025Subsec01}

Suppose the induced MDP is semicontinuous and has a feasible strategy with a finite value. Then in view of Theorem \ref{t1}, the relaxed problem (\ref{SashaLp02})-(\ref{e10}) can be rewritten as the following one:
\begin{eqnarray}\label{enum3}
\mbox{Minimize}&:&\int_{{\bf V} \times\textbf{B}}\bar{C}_0(x,b)\mu(dx\times db) ~\mbox{ over } \mu\in{\cal D}\\
  \mbox{subject to}&:& \int_{{\bf V} \times\textbf{B}}\bar{C}_j(x,b)\mu(dx\times db) - d_j\le 0,~~~j=1,2,\dots,J,\nonumber
\end{eqnarray}
where $\cal D$ is the set of measures  $\mu$ on ${\bf V}\times{\bf B}$ satisfying 
\begin{itemize}
\item[(a)]$\mu_{\bf V}\in {\cal M}^{\sigma<\infty}(\bf V)$, 
\item[(b)] equation (\ref{e9}), i.e., $\mu(dx\times \textbf{B})=\delta_{x_0}(dx)+\int_{{\bf V}\times \textbf{B}}Q(dx|y,b)\mu(dy\times db)$, and
\item[(c)]
\begin{equation}\label{e13p}
\int_{{\bf V}\times{\bf B}} \sum_{j=0}^J \bar C_j(x,b)\mu(dx\times db)<\infty.
\end{equation}
\end{itemize}
In this way, we may view  the relaxed problem (\ref{SashaLp02})-(\ref{e10}) as a convex program with $J$ constraint inequalities, because the set ${\cal D}$ is convex in the sense that if $c\in[0,1]$ and $\mu_1,\mu_2 \in{\cal D}$, then $c \mu_1+(1-c)\mu_2\in {\cal D}.$  It is also possible to embed ${\cal D}$ in a vector space, although this linear space is never needed in what follows.  We call this convex program the primal convex program. 

It is straightforward  to check that the primal convex program (\ref{enum3}) is equivalent to
\begin{equation}\label{enum4p}
\mbox{Minimize over }~\mu\in{\cal D}~:~~\sup_{\bar g\in\RR_+^J} L_1(\mu,\bar g),
\end{equation}
where  $\RR_+^J:=[0,\infty)^J$, and $$L_1(\mu,\bar g):=\int_{{\bf V} \times\textbf{B}}\bar{C}_0(x,b)\mu(dx\times db)+\sum_{j=1}^J g_j\left(\int_{{\bf V} \times\textbf{B}}\bar{C}_j(x,b)\mu(dx\times db)-d_j\right)$$ for $\mu\in{\cal D}$, $\bar g:=(g_1,g_2,\ldots, g_J)\in\RR_+^J$ is the Lagrangian, which is real-valued due to Condition (c) in the definition of ${\cal D}$. In greater detail, this follows from the following three observations. (a) for any $\mu\in {\cal D}$, if it is infeasible in  the primal convex program (\ref{enum3}), then  $\sup_{\bar g\in\RR_+^J} L_1(\mu,\bar g)=\infty$ for some $\bar{g}\in \RR_+^J$, meaning that in problem (\ref{enum4p}) it is sufficient to restrict to the set of $\mu\in {\cal D}$ feasible in the primal convex program (\ref{enum3}). (b) This set is non-empty, because  the optimal solution $\mu^*\in {\cal D}$ in the primal convex program (\ref{enum3}) exists by Theorem \ref{t1}(b). (c) If $\mu\in {\cal D}$ is feasible in the primal convex program (\ref{enum3}), then $\sup_{\bar g\in\RR_+^J} L_1(\mu,\bar g)=L_1(\mu,{\bf 0})
=\int_{{\bf V} \times\textbf{B}}\bar{C}_0(x,b)\mu(dx\times db),$ where ${\bf 0}=(0,0,\ldots,0)$. In particular, 
\begin{eqnarray}\label{PZ2025Eqn02}
val(\mbox{primal convex program (\ref{enum3})})=val(\mbox{problem (\ref{enum4p})}).
\end{eqnarray}
and both problems have the same optimal solutions.

The dual problem for the primal convex program  (\ref{enum4p}) reads
\begin{equation}\label{enum4}
\mbox{Maximize over }~ \bar g\in\RR_+^J~:~~\inf_{\mu\in{\cal D}} L_1(\mu,\bar g),
\end{equation}
It follows that $val(\mbox{problem (\ref{enum4})})\le val(\mbox{problem (\ref{enum4p})})$. There is a duality gap if the strict inequality holds. Next, in Theorem \ref{prop2}, we observe the absence of the duality gap, under the following Slater condition.

\begin{condition}[Slater condition]\label{con4} There exists $\mu\in{\cal D}$ such that
$$ \int_{{\bf V} \times\textbf{B}}\bar{C}_j(x,b)\mu(dx\times db) < d_j,~~~j=1,2,\dots,J.$$
\end{condition}
In view of Lemma \ref{l3}, given that the induced MDP is semicontinuous,  Condition \ref{con4} implies that the induced MDP has a feasible strategy with a finite value.

\begin{theorem}\label{prop2}
Suppose the induced MDP is semicontinuous and satisfies Condition \ref{con4} (i.e., the Slater condition).   Then the following assertions hold.
\begin{itemize}
\item[(a)]  $val(\mbox{primal convex program (\ref{enum4p}))}=\inf_{\mu\in{\cal D}}\sup_{\bar g\in\RR_+^J} L_1(\mu,\bar g)=\sup_{\bar g\in\RR_+^J} \inf_{\mu\in{\cal D}} L_1(\mu,\bar g)$,  there exists at least one $\bar g^*\in\RR_+^J$ solving the dual problem (\ref{enum4}), and the dual functional $h(\cdot)$ defined for each $\bar g\in\RR_+^J$ by $$h(\bar g):= \inf_{\mu\in{\cal D}} L_1(\mu,\bar g)$$ is concave.
\item[(b)] Consider a pair $(\mu^\ast,\bar{g}^\ast)\in {\cal D}\times \RR_+^J$. Then $\mu^*$ is an optimal  solution to the primal convex program  (\ref{enum4p}) and $\bar g^*$ is an optimal solution to the dual problem (\ref{enum4})
if and only if one of the following two equivalent statements holds.
\begin{itemize}
\item[(i)]  The pair $(\mu^*,\bar g^*)$ is a  saddle point of the Lagrangian:
\begin{eqnarray*}
L_1(\mu^*,\bar g)\le L_1(\mu^*,\bar g^*)\le L_1(\mu,\bar g^*),~\forall~ \mu\in{\cal D},~\bar g\in\RR_+^J.
\end{eqnarray*}
\item[(ii)]
The following relations hold for the pair $(\mu^*,\bar g^*)$:
\begin{eqnarray*}
&&\int_{{\bf V} \times\textbf{B}}\bar{C}_j(x,b)\mu^*(dx\times db) \le d_j,~j=1,2,\ldots, J;\\
&& L_1(\mu^*,\bar g^*)=\inf_{\mu\in{\cal D}} L_1(\mu,\bar g^*);\\
&& \sum_{j=1}^J g^*_j\left( \int_{{\bf V} \times\textbf{B}}\bar{C}_j(x,b)\mu^*(dx\times db)-d_j\right)=0.
\end{eqnarray*}
The last equality is known as the complementary slackness condition. 
\end{itemize}
\end{itemize}

\end{theorem}
\par\noindent\textit{Proof.} (a)  The absence of the duality gap and the solvability of the dual program (\ref{enum4}) follow from
Theorem 1, Section 8.6 of \cite{b42}, and the last assertion holds by Proposition 1, Section 8.6 of \cite{b42}. Note that in  the framework of Section 8.3 and Section 8.6 of \cite{b42}, ${\cal D}$ is embedded in a vector space. One can check that the reasoning therein for the results referred to in this proof holds if ${\cal D}$ is convex in the sense that if $c\in[0,1]$ and $\mu_1,\mu_2 \in{\cal D}$, then $c \mu_1+(1-c)\mu_2\in {\cal D}.$ 

(b) The first assertion in (b), i.e., the (primal and dual) optimality of $\mu^\ast$ and $\bar{g}^\ast$, is equivalent to  (b-i) according to (a) and Theorem 2 of \cite{rock}. Assume that   assertion (b-i) holds. Then, the first relation in  (b-ii) holds, because $\mu^\ast$ is optimal and in particular feasible in the primal convex program  (\ref{enum4p}). The second relation in  (b-ii) holds by (b-i).  The last relation in  (b-ii) holds by Theorem 1, Section 8.6 of \cite{b42}.  It remains to prove  that (b-ii) implies (b-i). Assume that (b-ii) holds. Then from $L_1(\mu^*,\bar g^*)=\inf_{\mu\in{\cal D}} L_1(\mu,\bar g^*)$ we see that  $ L_1(\mu^*,\bar g^*)\le L_1(\mu,\bar g^*)$ for all $\mu\in{\cal D}$.  From the complementary slackness condition, we see that  
 $L_1(\mu^*,\bar g^*)=\int_{{\bf V} \times\textbf{B}}\bar{C}_0(x,b)\mu^*(dx\times db)$. Now from 
 $\int_{{\bf V} \times\textbf{B}}\bar{C}_j(x,b)\mu^*(dx\times db) \le d_j,~j=1,2,\ldots, J$, we see that 
$
L_1(\mu^*,\bar g^*)=\int_{{\bf V} \times\textbf{B}}\bar{C}_0(x,b)\mu^*(dx\times db) \ge \int_{{\bf V}\times{\bf B}} \bar C_0(x,b)\mu^*(dx\times db)+\sum_{j=1}^J g_j\left(\int_{{\bf V}\times{\bf B}} \bar C_j(x,b)\mu^*(dx\times db)-d_j\right)=L_1(\mu^\ast,\bar g)
$
for all  $\bar g\in\RR_+^J$. Thus, \begin{eqnarray*}
L_1(\mu^*,\bar g)\le L_1(\mu^*,\bar g^*)\le L_1(\mu,\bar g^*),~\forall~ \mu\in{\cal D},~\bar g\in\RR_+^J,
\end{eqnarray*}
i.e.,  (b-ii) implies (b-i). $\hfill\Box$

 \section{Main results}\label{PZ2025Section02}
\subsection{Second convex program and duality results}\label{PZ2025Subsec02}

\begin{condition}\label{con5}\begin{itemize}
\item[(a)] $C^I_0(x,a)\ge \delta>0$ for all $x\in{\bf X},~a\in{\bf A}$.
\item[(b)]  ${\cal C}=\max_{j=0,1,\ldots,J}\sup_{x\in{\bf X}}C^g_j(x)+\max_{j=0,1,\ldots,J}\sup_{(x,a)\in{\bf X}\times{\bf A}}C^I_j(x,a)<\infty$.
\end{itemize}
\end{condition}

Suppose the induced MDP is semicontinuous. Under Condition \ref{con5}(a), for each $x\in{\bf V}^c\setminus\{\Delta\}$, by (\ref{e5}), necessarily  $f^*(x)=(\infty,\hat a)$ for some $\hat a\in{\bf A}$. Hence, under Condition \ref{con5}, for each $x\in{\bf V}^c\setminus\{\Delta\}$ $C^g_j(\phi(x,t))= 0$ for almost all $t\ge 0$ with respect to the Lebesgue measure for all $j=0,1,\ldots, J$. One can also say that, for each given $x\in {\bf X},$ $x\in{\bf V}^c$ if and only if $C^g_j(\phi(x,t))\equiv 0$ a.e. with respect to the Lebesgue measure for all $j=0,1,\ldots,J$, and $Q(\{\Delta\}|x,f^*(x))=1$.

\begin{lemma}\label{l4}
Suppose Condition \ref{con5}(a) is satisfied. For each given strategy $\pi$, if $\mu^\pi({\bf X}\times{\bf B})=\infty$, then ${\cal V}(x_0,\pi,\bar C_0)=\infty$.
\end{lemma}

Suppose that the induced MDP is semicontinuous. Recall from Theorem \ref{t1} that the impulse control problem (\ref{PZZeqn02}) is equivalent to the relaxed problem (\ref{SashaLp02})-(\ref{e10}). This problem has a finite optimal value, and is solvable. In view of Lemma \ref{l4}, if in addition, Condition \ref{con5}(a) is also satisfied, then problem  (\ref{SashaLp02})-(\ref{e10}) is equivalent to the following convex program: 
\begin{eqnarray}\label{e21}
\mbox{Minimize}&:&\int_{{\bf V} \times\textbf{B}}\bar{C}_0(x,b)\mu(dx\times db)
 \mbox{ over $\mu \in {\cal M}^{<\infty}({\bf V}\times {\bf B})$} \\
 \mbox{subject to}&:& \mu(\Gamma\times \textbf{B})=\delta_{x_0}(\Gamma)+\int_{{\bf V}\times \textbf{B}}Q(\Gamma|y,b)\mu(dy\times db),~\Gamma\in{\cal B}({\bf V})\label{e22}\\
\mbox{and}&&\int_{{\bf V} \times\textbf{B}}\bar{C}_j(x,b)\mu(dx\times db)\le d_j,~j=1,2,\dots,J, \label{e23}
\end{eqnarray}
where ${\cal M}^{<\infty}({\bf V}\times {\bf B})$ is the set of finite measures on $({\bf V}\times {\bf B},{\cal B}({\bf V}\times {\bf B}))$. We shall refer to problem  (\ref{e21})-(\ref{e23}) as the second (primal) convex program with one measure-valued constraint equality and $J$ constraint inequalities to distinguish it from the primal convex program (\ref{enum3}) with $J$ constraint inequalities. Note that (\ref{e21})-(\ref{e23}) is the linear program in the linear space of finite signed measures on $({\bf V}\times{\bf B})$.

Suppose that the induced MDP is semicontinuous  and has a feasible strategy with a finite value, and Condition \ref{con5}(a) is satisfied. We shall formulate a dual problem for the second convex program (\ref{e21})-(\ref{e23}) and show the absence of the duality gap. Corresponding to the constraint equality in problem  (\ref{e21})-(\ref{e23}), let us introduce $\mathbb{B}(\bf V)$, the set of bounded measurable functions on ${\bf V}.$ Then we introduce the Lagrangian $L_2(\cdot,\cdot,\cdot)$ defined for each $(\mu,W(\cdot),\bar g)\in {\cal M}^{<\infty}({\bf V}\times {\bf B})\times \mathbb{B}({\bf V})\times \RR_+^J$ by
\begin{eqnarray*}
L_2(\mu,W,\bar g)&:=& \int_{{\bf V}\times{\bf B}} \bar C_0(x,b)\mu(dx\times db)\\
&&+\int_{\bf V} W(x)\left[\delta_{x_0}(dx)+\int_{\bf V} Q(dx|y,b)\mu(dy\times db)-\mu(dx\times{\bf B})\right]\\
&&+\sum_{j=1}^J g_j\left(\int_{{\bf V}\times{\bf B}}\bar C_j(x,b)\mu(dx\times db)-d_j\right).
\end{eqnarray*}
Since $\bar C_j(\cdot,\cdot)$ and $g_j$ are non-negative, the Lagrangian $L_2(\cdot,\cdot,\cdot)$ is well defined and the terms therein can be arbitrarily rearranged, as   $\int_{{\bf V}\times{\bf B}}\bar C_j(x,b)\mu(dx\times db)$ are the only possibly infinite terms, which are non-negative and with non-negative coefficients. Here $0\times \infty:=0$ is accepted.

It is straightforward to check that the second primal convex (in fact, linear) program (\ref{e21})-(\ref{e23}) is equivalent to
\begin{equation}\label{e24}
\mbox{Minimize over }~\mu\in{\cal M}^{<\infty}({\bf V}\times{\bf B}):~~\sup_{(W(\cdot),\bar g)\in{\mathbb{B}({\bf V})}\times\RR_+^J} L_2(\mu,W,\bar g).
\end{equation}
In fact, it follows from the observation formulated in the next remark.

\begin{remark}\label{r2}
If equality (\ref{e22}) is violated by some $\mu\in{\cal M}^{<\infty}({\bf V}\times{\bf B})$, then there exist $\varepsilon>0$ and $\Gamma\in{\cal B}({\bf V})$ such that
$$\delta_{x_0}(\Gamma)+\int_{\bf  V} Q(\Gamma|y,b)\mu(dy\times db)-\mu(\Gamma\times{\bf B})\ge\varepsilon~(\mbox{or } \le-\varepsilon).$$
Now for $W_N(x)=N\II\{x\in\Gamma\}$ (resp., for $W_N(x)=-N\II\{x\in\Gamma\}$),  
$\lim_{N\to\infty} L_2(\mu,W_N,\bar g)=+\infty$. The similar reasoning applies to the constraints (\ref{e23}). The conclusion is that in problem (\ref{e24}), one can be restricted to the set of measures $\mu$ that are feasible in the second primal convex program (\ref{e21})-(\ref{e23}),  and  problem (\ref{e24}) restricted on this set becomes the same as the second primal convex program (\ref{e21})-(\ref{e23}).
\end{remark}

Now, we formulate the dual  program for the second primal convex program (\ref{e21})-(\ref{e23}) as 
\begin{equation}\label{e25}
\mbox{Maximize over }~(W(\cdot),\bar g)\in \mathbb{B}({\bf V})\times\RR_+^J:~\inf_{\mu\in {\cal M}^{<\infty}({\bf V}\times {\bf B})} L_2(\mu,W,\bar g).
\end{equation}

One can solve the dual  program (\ref{e25}) in two steps: firstly, for each fixed $\bar g\in\RR_+^J$, solve problem
\begin{equation}\label{e26}
\mbox{Maximize over $W(\cdot)\in\mathbb{B}({\bf V})$: } \inf_{\mu\in {\cal M}^{<\infty}({\bf V}\times{\bf B})} L_2(\mu,W,\bar g);
\end{equation}
and after that  solve problem 
\begin{eqnarray*}
\mbox{Maximize over $\bar{g}\in \RR_+^J$: } \sup_{W(\cdot)\in \mathbb{B}({\bf V})} \inf_{\mu\in {\cal M}^{<\infty}({\bf V}\times{\bf B})} L_2(\mu,W,\bar g).
\end{eqnarray*}
Up to the additive  constant $-\sum_{j=1}^J g_jd_j$, problem (\ref{e26}) can be rewritten as follows:
\begin{eqnarray}\label{e27}
\mbox{Maximize}&:& W(x_0)
 \mbox{ over } W(\cdot)\in\mathbb{B}({\bf V}) \\
 \mbox{subject to}&:&
 \bar C_0(x,b)+\sum_{j=1}^J g_j\bar C_j(x,b)+\int_{\bf V} W(y) Q(dy|x,b)-W(x)\ge 0,\label{e28}\\
 &&~~~~~~~~~~~(x,b)\in{\bf V}\times{\bf B}.\nonumber
\end{eqnarray}
Indeed, if inequality (\ref{e28}) is violated by some $W(\cdot)\in\mathbb{B}({\bf V})$, then, for some $\varepsilon>0$, $ \bar C_0(\hat{x},\hat{b})+\sum_{j=1}^J g_j\bar C_j(\hat{x},\hat{b})+\int_{\bf V} W(y) Q(dy|\hat{x},\hat{b})-W(\hat{x})<-\varepsilon$  for some $(\hat x,\hat b)\in{\bf V}\times{\bf B}$. Now, for the sequence of finite measures defined for each $N\in\{1,2,\dots\}$ by $\mu_N(dx\times db)=N\delta_{(\hat x,\hat b)}(dx\times db)\in {\cal M}^{<\infty}({\bf V}\times {\bf B})$,   $\lim_{N\to\infty}L_2(\mu_N,W,\bar g)=-\infty$, and thus $\inf_{\mu\in {\cal M}^{<\infty}({\bf V}\times {\bf B})} L_2(\mu,W,\bar g)=-\infty.$
On the other hand, for each  $W(\cdot)\in\mathbb{B}({\bf V})$ satisfying (\ref{e28}), $\inf_{\mu\in {\cal M}^{<\infty}({\bf V}\times {\bf B})} L_2(\mu,W,\bar g)$  is provided by the null measure $\mu\equiv 0$ and equals $W(x_0)-\sum_{j=1}^J g_jd_j$.

It is observed in Theorem \ref{t5} below that for each $\bar g\in \RR_+^J$,  problem (\ref{e27})-(\ref{e28}) is solvable. In fact, under the conditions imposed in Theorem \ref{t5}, its optimal solution is given by the function $W^\ast_{\bar g}(\cdot)$ 
defined by
\begin{eqnarray}\label{e31}
W^*_{\bar g}(x)&:=&\inf_\pi\EE^\pi_x\left[\sum_{i=1}^\infty\left(\bar C_0(X_{i-1},B_i)+\sum_{j=1}^J g_j\bar C_j(X_{i-1},B_i)\right)\right]\\
&=&\inf_\pi\left\{{\cal V}(x,\pi,\bar{C}_0)+\sum_{j=1}^J g_j{\cal V}(x,\pi,\bar{C}_j)\right\},
~~~x\in{\bf V}.\nonumber
\end{eqnarray}
Thus, $W^*_{\bar g}(\cdot)$ is the restriction on ${\bf V}$ of the Bellman function  ${\cal V}^\ast(\cdot, \bar{C}_0+\sum_{j=1}^J g_j \bar{C}_j)$ defined for each $x\in {\bf X}_\Delta$ by  ${\cal V}^\ast(x, \bar{C}_0+\sum_{j=1}^J g_j \bar{C}_j)(x):= \inf_\pi {\cal V}(x,\pi, \bar{C}_0+\sum_{j=1}^J g_j \bar{C}_j)$. We also call $W^*_{\bar g}(\cdot)$ the Bellman function because ${\cal V}^\ast(x, \bar{C}_0+\sum_{j=1}^J g_j \bar{C}_j)=0$ for all $x\in {\bf V}^c.$
The next lemma shows that  $W^\ast_{\bar{g}}(\cdot)$ is a feasible solution in  problem (\ref{e27})-(\ref{e28}) and provides a useful characterization of the function $W^\ast_{\bar g}(\cdot)$.

\begin{lemma}\label{l7}
Suppose the induced MDP is semicontinuous  and Condition \ref{con5} is  satisfied. Let $\bar g\in\RR_+^J$ is arbitrarily fixed. Then the following assertions hold.
\begin{itemize}
\item[(a)] $W^\ast_{\bar{g}}(\cdot)$ is bounded, non-negative and lower semicontinuous on ${\bf V}$, and $W^\ast_{\bar{g}}(\cdot)$ is a feasible solution in  problem (\ref{e27})-(\ref{e28}).
\item[(b)] $W^\ast_{\bar{g}}(\cdot)$ is the unique lower semicontinuous bounded function on ${\bf V}$ satisfying the following equation:
\begin{eqnarray*}
W(x) &=& \inf_{b\in{\bf B}}\left\{\bar C_0(x,b)+\sum_{j=1}^J g_j\bar C_j(x,b)+\int_{\bf V} W(y) Q(dy|x,b)\right\},~\forall~x\in{\bf V}
\end{eqnarray*}
Moreover,
if a function $W(\cdot)\in\mathbb{B}({\bf V})$ satisfies
\begin{eqnarray}\label{enn}
W(x) &=& \inf_{b\in{\bf B}}\left\{\bar C_0(x,b)+\sum_{j=1}^J g_j\bar C_j(x,b)+\int_{\bf V} W(y) Q(dy|x,b)\right\}\\
&=& \bar C_0(x,\hat{f}(x))+\sum_{j=1}^J g_j\bar C_j(x,\hat{f}(x))+\int_{\bf V} W(y) Q(dy|x,\hat{f}(x)),~~~x\in{\bf V}\nonumber
\end{eqnarray}
for a measurable mapping $\hat{f}(\cdot):~{\bf V}\to{\bf B}$,
then the deterministic stationary strategy $\hat{f}$ (supplemented by $\hat{f}(x):=f^*(x)$ for $x\in{\bf V}^c$, see Definition \ref{d3}) is uniformly optimal in the sense that
$${\cal V}(x,\hat{f},\bar{C}_0)+\sum_{j=1}^J g_j{\cal V}(x,\hat{f},\bar{C}_j)={\cal V}^\ast(x, \bar{C}_0+\sum_{j=1}^J g_j \bar{C}_j),~~~x\in{\bf X}_\Delta.$$ 
\end{itemize}
\end{lemma}

\begin{theorem}\label{t5}
Suppose the induced MDP is semicontinuous, and Conditions \ref{con4} and \ref{con5} are satisfied. Then the following assertions hold.
\begin{itemize}
\item[(a)] For a fixed $\bar g\in\RR_+^J$, the Bellman function $W^\ast_{\bar{g}}(\cdot)$ defined by (\ref{e31})
is bounded, non-negative and lower semicontinuous, and
solves problem (\ref{e27})-(\ref{e28}) with $$val(\mbox{problem (\ref{e27})-(\ref{e28})})=W^*_{\bar g}(x_0).$$
\item[(b)]  $W^*_{\bar g}(x_0)-\sum_{j=1}^J g_jd_j=h(\bar g)$, where $h(\cdot)$ was defined in Theorem \ref{prop2}(a). Thus, a maximizer $\bar g^*\in \RR_+^J$ for
$$\max_{\bar g\in\RR_+^J}\left\{\max_{W(\cdot)\in \mathbb{B}({\bf V})} \inf_{\mu\in{\cal M}^{<\infty}({\bf V}\times{\bf B})} L_2(\mu,W,\bar g)\right\}=\max_{\bar g\in\RR_+^J} h(\bar g)$$
exists.
\item[(c)] The pair $(W^*_{\bar g^*},\bar g^*)$ is a solution to the dual program (\ref{e25}).
\item[(d)] There is no duality gap:
$$\inf_{\mu\in{{\cal M}^{<\infty}({\bf V}\times {\bf B})}}\sup_{(W(\cdot),\bar g)\in{\mathbb{B}({\bf V})}\times\RR_+^J} L_2(\mu,W,\bar g)=\sup_{(W(\cdot),\bar g)\in{\mathbb{B}({\bf V})}\times\RR_+^J}\inf_{\mu\in{{\cal M}^{<\infty}({\bf V}\times {\bf B})}}L_2(\mu,W,\bar g).$$
\end{itemize}
\end{theorem}

In view of (\ref{e31}) and the definition of $h(\bar{g})=\inf_{\mu\in {\cal D}}L_1(\mu,\bar{g})$ given in Theorem \ref{prop2}(a), under the conditions imposed in Theorem \ref{t5}, part (b) therein asserts that 
\begin{eqnarray*}
\inf_\pi\EE^\pi_x\left[\sum_{i=1}^\infty\left(\bar C_0(X_{i-1},B_i)+\sum_{j=1}^J g_j\bar C_j(X_{i-1},B_i)\right)\right]-\sum_{j=1}^J g_jd_j=\inf_{\mu\in {\cal D}}L_1(\mu,\bar{g})
\end{eqnarray*}
for any $\bar{g}\in \RR_+^J$. This equality is attractive from the practical point of view, because it shows that the value of the dual functional $h(\bar{g})$ can be calculated using the dynamic programming approach, and this is further elaborated on in the next subsection. On the other hand, here we remark that Condition \ref{con5} is important for the last equality. 

\subsection{Procedures for finding $(J+1)$-mixed optimal strategies in the original impulse control problem}\label{sec6}
Throughout this subsection, we suppose that the induced MDP is semicontinuous, and Conditions \ref{con4} and \ref{con5} are satisfied, so that all the above results are applicable. In view of Theorem \ref{t5}, we now present a general procedure for solving the primal convex problem  (\ref{enum3})  and the relaxed problem (\ref{SashaLp02})-(\ref{e10}), leading to optimal $(J+1)$-mixed strategies, whose definition is given above Proposition \ref{t2}, for the original impulse control problem (\ref{PZZeqn02}). The main result in this subsection is the justification of this procedure, see Theorem \ref{l8} below. An application of this procedure is illustrated in the next section.

\par\noindent\textit{Formulation of the procedure.}
\par\noindent\textit{Step 1.} For each fixed $\bar g\in\RR_+^J$, compute the dual functional  $h(\bar g)=\inf_{\mu\in{\cal D}}L_1(\mu,\bar g)$. 
According to Theorem \ref{t5}(b), $W^*_{\bar g}(x_0)-\sum_{j=1}^J g_jd_j=h(\bar g)$. Thus, it suffices to compute the Bellman function $W^*_{\bar g}(\cdot)$, which is bounded. To this end, since the induced MDP is positive and semicontinuous, one may apply the successive approximations, see Corollary 9.17.2 of \cite{Bertsekas:1978}.  

\par\noindent\textit{Step 2.} Find any $\bar g^*$ providing $\max_{\bar g\in\RR_+^J} h(\bar g)$.  By Theorem \ref{prop2}(a) and Theorem \ref{t5}(b),  such $\bar{g}^\ast$ exists and is optimal in problem (\ref{enum4}), and the function $h(\cdot)$ is concave.  In this step, one can use any standard method of maximizing  concave functions on $\RR_+^J$.

\par\noindent\textit{Step 3.} For $\bar g^*$ obtained in Step 2,  calculate the Bellman function $W^*_{\bar g^*}(\cdot)$ as in Step 1. 
Let  ${\bf F}$ denote the set of  deterministic stationary strategies $f$ providing the minimum in the Bellman equation:
\begin{eqnarray*}
W^*_{\bar g^*}(x) &=& \inf_{b\in{\bf B}}\left\{\bar C_0(x,b)+\sum_{j=1}^J g^*_j\bar C_j(x,b)+\int_{\bf V} W^*_{\bar g^*}(y)Q(dy|x,b)\right\}\\
&=& \bar C_0(x,f(x))+\sum_{j=1}^J g^*_j\bar C_j(x,f(x))+\int_{\bf V} W^*_{\bar g^*}(y)Q(dy|x,f(x)),~x\in{\bf V}
\end{eqnarray*}
and satisfying $f(x)=f^*(x)$ for $x\in{\bf V}^c.$ Note that the expression in the parentheses in the first line of the above equality is a lower semicontinuous function and $\bf B$ is compact, so that ${\bf F}\ne \emptyset,$  see Proposition 7.33 of  \cite{Bertsekas:1978} and Theorem 3.3 of \cite{Feinberg2013JMAA}. Moreover, 
according to Proposition 9.12 of \cite{Bertsekas:1978}, each $f\in {\bf F}$ is uniformly optimal in the sense that ${\cal V}(x,f,\bar{C}_0+\sum_{j=1}^J g_j^\ast \bar{C}_j)=\inf_\pi{\cal V}(x,\pi,\bar{C}_0+\sum_{j=1}^J g_j^\ast \bar{C}_j)=W^*_{\bar g^*}(x)$ for all $x\in {\bf X}_\Delta.$

\par\noindent\textit{Step 4.} Find constants $\gamma_l\in (0,1]$ and deterministic stationary strategies $f^l\in {\bf F}$ (not necessarily all distinct) for $l=1,2,\ldots,J+1$ satisfying 
\begin{itemize}
\item[(a)] $\sum_{l=1}^{J+1}\gamma_l=1$;
\item[(b)] for the convex combination $\mu^*=\sum_{l=1}^{J+1} \gamma_l\mu^{f^l}|_{{\bf V}\times {\bf B}}$ of the restrictions of the occupation measures of $f^l$, it holds for each $j=1,2,\ldots,J$ that $\int_{{\bf V}\times{\bf B}} \bar C_j(x,b)\mu^*(dx\times db)\le d_j$, whereas  $\int_{{\bf V}\times{\bf B}} \bar C_j(x,b)\mu^*(dx\times db)= d_j$ if $g_j^*>0$.
\end{itemize}
Provided that we can find such constants $\gamma_l$ and strategies $f^l,$ $\mu^*$  would solve the primal convex program  (\ref{enum3})  as well as the relaxed problem (\ref{SashaLp02})-(\ref{e10}). This can be seen as follows. From the proof of  Theorem \ref{t5}(a), we see for each $f\in{\bf F}$ that  $\mu^f({\bf V}\times{\bf B})<\infty$, and from  the inequality $\int_{{\bf V}\times{\bf B}}\bar C_j(x,b)\mu^*(dx\times db)\le d_j$ (and keeping in mind that $\gamma_l>0$ for all $l=1,2,\ldots,J+1$) we deduce that, for all $j=1,2,\ldots,J$,
$${\cal V}(x_0,f^l,\bar{C}_j)=\int_{{\bf V}\times{\bf B}} \bar C_j(x,b)\mu^{f^l}(dx\times db)<\infty.$$
Therefore, $\mu^{f^l}\in{\cal D}$ for all $l=1,2,\ldots, J+1$ and $\mu^*\in{\cal D}$, too. Furthermore, from the proof of Theorem \ref{t5}(b), we see that
$$L_1(\mu^{f^l}|_{{\bf V}\times{\bf B}},\bar g^*)=\inf_{\mu\in{\cal D}} L_1(\mu,\bar g^*)=W^*_{\bar g^*}(x_0)-\sum_{j=1}^J g^*_j d_j.$$ 
It follows that
$$L_1(\mu^*,\bar g^*)=\sum_{l=1}^{J+1} \gamma_l L_1(\mu^{f^l},\bar g^*)=\inf_{\mu\in{\cal D}} L_1(\mu,\bar g^*).$$ 
Thus, the pair $(\mu^\ast,\bar{g}^\ast)\in {\cal D}\times \RR_+^J$ satisfies the conditions in Theorem \ref{prop2}(b-ii), and the optimality of $\mu^\ast$ now follows from Theorem \ref{prop2}(b). 
Moreover, $\mu^\ast=\mu^{\pi^\ast}|_{{\bf V}\times {\bf B}}$ for some strategy $\pi^\ast$. Such a strategy is called  a mixture of the deterministic stationary strategies $f^l$. Its existence follows from the next fact, which follows from the convexity of the space of all strategic measures, see Chapter 5, Section 5 of \cite{dyn}; see also Propositions 1 and 3 of \cite{sicon24}. 
\begin{proposition}\label{l5}
The set of all occupation measures in the induced MDP is convex, in the sense that if $c\in[0,1]$ and $\mu_1$ and $\mu_2$ are occupation measures in the induced MDP, then $c \mu_1+(1-c)\mu_2$ is the occupation measure of some strategy $\pi.$
\end{proposition}

The next result shows that one can always find such constants $\gamma_l$ and strategies $f^l$ as in Step 4.

\begin{theorem}\label{l8}
Suppose the induced MDP is semicontinuous, Conditions  \ref{con4} and \ref{con5} are satisfied, and ${\bf V}={\bf X}$.  Then
there exists a  mixture $\pi^\ast\in \Pi$ of $J+1$ deterministic stationary strategies $f^l\in{\bf F}$ with coefficients $\{\gamma_l\}_{l=1}^{J+1}\subseteq(0,1]$, such that $\pi^*$ is optimal  in the original impulse control problem (\ref{PZZeqn02}) and $\mu^{\pi^\ast}|_{{\bf V}\times {\bf B}}$ solves the primal convex program  (\ref{enum3})  and the relaxed problem (\ref{SashaLp02})-(\ref{e10}).
\end{theorem}
Consider $\pi^\ast$ as in Theorem \ref{l8}. Let $\mu^\ast:=\mu^{\pi^\ast}|_{{\bf V}\times {\bf B}}$. Since $\mu^\ast$ solves the primal convex program  (\ref{enum3}), for any $\bar{g}^\ast\in \RR_+^J$ providing the maximum of the dual functional $h(\bar{g})$, we see from Theorem \ref{prop2}(b) that condition (b) in Step 4 of the general procedure is satisfied by $\mu^\ast$, and   there exist constants $\{\gamma_l\}_{l=1}^{J+1}\subseteq(0,1]$ and strategies $\{f^l\}_{l=1}^{J+1}\subset{\bf F}$  satisfying the both conditions (a) and (b) in Step 4.

The proof of Theorem \ref{l8} is based on the folowing statements.

\begin{proposition}\label{t2}
Suppose the induced MDP is semicontinuous,   there exists a feasible strategy with a finite value, Condition  \ref{con5}(a) is satisfied, and ${\bf V}={\bf X}$. Then there exists an optimal strategy $\pi^*$ in the original impulse control problem (\ref{PZZeqn02})  in the form of a mixture of $J+1$ deterministic stationary strategies $f^1,f^2,\ldots,f^{J+1}$.
\end{proposition}

For the proof see Theorem 2 of \cite{sicon24}. Note that the strategies $f^l$ here may be not uniformly optimal.
 
Consequently, the key component in the proof of Theorem \ref{l8} is to find the suitable versions in ${\bf F}$  of $f^l$ in Proposition \ref{t2}. This is based on the following result for generic MDPs,  which may be of independent interest, where ${\bf V}$ is generic, not necessarily generated by the primitives of the impulse control model.

\begin{lemma}\label{l6} 
Consider a generic MDP $\{{\bf Y},{\bf B},Q,C\}$, where ${\bf Y}={\bf V}\cup\{\Delta\}$ with $\Delta\notin{\bf V}$ being the  isolated cemetery. Suppose the single cost function $C(\cdot,\cdot)$ is $[0,\infty]$-valued, $C(\Delta,b)\equiv 0$, and this MDP is semicontinuous. Let the Bellman function  
$${\cal V}^\ast(y,C):= \inf_\pi  \EE^\pi_y\left[\sum_{i=1}^\infty C(Y_{i-1},B_i)\right]$$
be defined for each $y\in {\bf Y}$. Clearly, ${\cal V}^*(\Delta,C)=0$. Let the initial state $y_0\in{\bf V}$ be fixed, and $f$ be an (optimal) deterministic stationary strategy such that
$$\EE^f_{y_0}\left[\sum_{i=1}^\infty C(Y_{i-1},f(Y_{i-1}))\right]=\int_{\bf V} C(y,f(y))\mu^f(dy\times{\bf B})={\cal V}^\ast(y_0,C),$$
where $\mu^f$ is the occupation measure. When $y=\Delta$, we fix $f(\Delta)=f^*(\Delta)$ arbitrarily.  Assume that ${\cal V}^\ast(y_0,C)<\infty.$
Then $\mu^f=\mu^{\tilde f}$ on ${\bf V}\times{\bf B}$ for some deterministic stationary strategy $\tilde f$, which is  uniformly optimal, i.e.,
\begin{equation}\label{e11}
\EE^{\tilde f}_y\left[\sum_{i=1}^\infty C(Y_{i-1},\tilde f(Y_{i-1}))\right]={\cal V}^\ast(y,C),~~~y\in{\bf V}.
\end{equation}
Similarly to the above, $\tilde f(\Delta):=f^*(\Delta)$ without loss of generality.
\end{lemma}

\section{Illustrative example: fluid queueing system}\label{sec7}
In this section, we illustrate the main results in the last section with an example. 
\begin{example}[fluid queueing system]\label{PZ2025Example02}
Consider the impulse control problem with the following system primitives. Let the state space be ${\bf X}=[0,\infty)$, and the action space ${\bf A}=\{1\}$ be a singleton, so that we omit the action $a\in {\bf A}$ in the description of the model below. The flow is described by $\phi(x,t)=x+t$ for all $x\in {\bf X}$ and $t\in[0,\infty).$ The new state upon applying an impulse is $l(x)\equiv 0$. Let $\alpha>0$, $x_0=0$. We consider the case of a single constraint, i.e., $J=1$, so that the index $j=1$ will be omitted below. Let the constraint constant be  $d>0$.  The impulse cost functions and running cost rates are given by 
$$C^I_0(x)\equiv K>0,~C^g_0(x)\equiv 0,~C^I_1(x)\equiv 0,~C^g_1(x)\equiv h>0,$$ where $h\in(0,\infty)$ and $K\in(0,\infty)$ are fixed constants. 
\end{example}

\par\noindent\textbf{Interpretation of the example.} This impulse control problem formulated in Example \ref{PZ2025Example02} can be given the following interpretation. Consider an infinite buffer: small packets of information arriving at a high speed are served with the high intensity. In such a case, this queueing system is known to be well approximated by the so called fluid model
$\dot x=\lambda-\mu:$
see e.g., \cite{mm}. Here $x\ge 0$ is the buffer filling size measured, e.g., in Mb and $\lambda$ (or $\mu$) is the arrival (resp., service) rate measured in Mb/time unit. We assume that $\lambda>\mu$ and choose the time scale such that $\lambda-\mu=1$. Thus, the flow is described by $\phi(x,t)=x+t$ and the state space is ${\bf X}=[0,\infty)$.
At any time moment the decision maker can use an additional very powerful server; as the result, the current state instantly changes to zero. Thus, we put ${\bf A}=\{1\}$. The result of applying the impulse is described by the constant function $l(x)\equiv 0$. The price for one impulse is $K>0$, and the holding cost is $h>0$ per unit per time unit. $\alpha>0$ is the discount factor as usual. We will fix the initial state $x_0=0$ and minimize the total discounted cost of impulses under the constraint $d>0$ on the total discounted holding cost.

\par\noindent\textbf{Induced MDP.}
One can compute the cost functions and write down the transition kernel in the induced MDP as follows: 
\begin{eqnarray*}
\bar C_0(x,\theta) &=& \left\{\begin{array}{ll} Ke^{-\alpha\theta} &\mbox{ if } \theta<\infty;\\
0, &\mbox{ if } \theta=\infty,\end{array}\right.\\
\bar C_1(x,\theta) &=& \left\{\begin{array}{ll} h\int_0^\theta e^{-\alpha t}(x+t)dt
=\frac{hx}{\alpha}\left(1-e^{-\alpha\theta}\right)+h\left[\frac{1}{\alpha^2}-\frac{1}{\alpha^2} e^{-\alpha\theta}-\frac{\theta}{\alpha} e^{-\alpha\theta}\right], & \mbox{ if } \theta<\infty;\\
\frac{hx}{\alpha}+\frac{h}{\alpha^2}, & \mbox{ if } \theta=\infty,\end{array}\right.\\
Q(dy|x,\theta)&=& \left\{\begin{array}{ll} e^{-\alpha\theta}\delta_0(dy)+(1-e^{-\alpha\theta})\delta_\Delta (dy) & \mbox{ if } \theta<\infty,~x\ne\Delta;\\
\delta_\Delta(dy) & \mbox{ otherwise}.\end{array}\right.
\end{eqnarray*}
Since there is a single impulse action, which we omit for brevity, a strategy is identified with the sequence of stochastic kernels $\pi_i(d\theta|x_0,\theta_1,x_1,\dots,x_{i-1})$ on $[0,\infty]$.

Condition \ref{ConstrainedPPZcondition01} is obviously satisfied, and thus the induced MDP is semicontinuous by Proposition \ref{pr1}. Condition \ref{con5} is obviously satisfied, too.  According to the comments below Condition \ref{con5}, ${\bf V}={\bf X}$. In order to check Condition \ref{con4}, consider the deterministic stationary strategy $$\hat f(x)=\left\{\begin{array}{ll} \hat\theta-x, & \mbox{ if } x\le\hat\theta;\\0, & \mbox{ if } x>\hat\theta,\end{array}\right.$$ i.e., apply the impulse as soon as $x\ge \hat\theta$, where $\hat\theta>0$ is fixed. Then, during one cycle formed by two consecutive visits to state $0$, the holding cost equals $\bar C_1(0,\hat\theta)$, and
$${\cal V}(0,\hat f,\bar{C}_1) = \sum_{i=0}^\infty e^{-\alpha\hat\theta i} \bar C_1(0,\hat\theta)=\frac{h\left[\frac{1}{\alpha^2}-\frac{1}{\alpha^2} e^{-\alpha\hat{\theta}}-\frac{\hat\theta}{\alpha} e^{-\alpha\hat{\theta}}\right]}{1-e^{-\alpha\hat\theta}}.$$
The above quantity depends on $\hat{\theta}$, and one can compute its limit
$$\lim_{\hat{\theta}\downarrow 0}\frac{h\left[\frac{1}{\alpha^2}-\frac{1}{\alpha^2} e^{-\alpha\hat{\theta}}-\frac{\theta}{\alpha} e^{-\alpha\hat{\theta}}\right]}{1-e^{-\alpha\hat\theta}}=\lim_{\hat\theta\downarrow 0} h\frac{\frac{1}{\alpha}e^{-\alpha\hat\theta}-\frac{1}{\alpha}e^{-\alpha\hat\theta}+\theta e^{-\alpha\hat\theta}}{\alpha e^{-\alpha\hat\theta}}=0.$$
Therefore, for the constraint constant $d>0$, when choosing small enough $\hat\theta>0$, one can satisfy the constraint ${\cal V}(0,\hat f,\bar{C}_1)<d$. In the meanwhile, similarly to the above,
$${\cal V}(0,\hat f,\bar{C}_0)=\frac{Ke^{-\alpha\hat\theta}}{1- e^{-\alpha\hat\theta}}<\infty.$$ In view of Lemma \ref{l2}, we see that Condition \ref{con4} (i.e., the Slater condition) holds for $\mu=\mu^{\hat f}|_{{\bf V}\times {\bf B}}$ with a sufficiently small $\hat\theta>0$; $\mu(dx\times{\bf B})=\delta_0(dx)\frac{1}{1- e^{-\alpha\hat\theta}}\in{\cal D}$ is a finite measure.

\par\noindent\textbf{Solving the impulse control problem in Example \ref{PZ2025Example02}.}  We now solve the impulse control problem in Example \ref{PZ2025Example02} by following the procedures in Subsection \ref{sec6}. We have seen that all the conditions in that subsection are satisfied in this example. Since $J=1$, we omit the bar in the Lagrange multiplier and write just $g$ and $g^*$.

\textbf{Step 1.} We find $W^\ast_g(\cdot)$ for $g\in[0,\infty).$ 

First, let $g>0.$ In view of  Lemma \ref{l7}, let us solve for the unique nonnegative bounded lower semicontinuous solution to the following equation: 
\begin{eqnarray}
W(x)& = & \inf_{\theta\in[0,\infty]} \left\{\bar C_0(x,\theta)+g\bar C_1(x,\theta)+\int_{\bf X} W(y)Q(dy|x,\theta)\right\}\nonumber\\
&=&\inf_{\theta\in[0,\infty]}\tilde{G}(x,\theta),~~~~~x\in{\bf V}, \label{eqn6}
\end{eqnarray}
where $\tilde{G}(x,\theta)=K e^{-\alpha \theta}+\frac{ghx}{\alpha}(1-e^{-\alpha \theta})+gh\left\{\frac{1}{\alpha^2}-\frac{1}{\alpha^2}e^{-\alpha\theta}-\frac{\theta}{\alpha}e^{-\alpha \theta}\right\}+e^{-\alpha\theta}W(0).$

When $x=0$, we have 
\begin{eqnarray*} 
&&W(0) =\inf_{\theta\in[0,\infty]}\left\{K e^{-\alpha \theta}+gh\left\{\frac{1}{\alpha^2}-\frac{1}{\alpha^2}e^{-\alpha\theta}-\frac{\theta}{\alpha}e^{-\alpha \theta}\right\}+e^{-\alpha\theta}W(0)\right\}.
\end{eqnarray*}
Let us denote $G(\theta)=K e^{-\alpha \theta}+gh\left\{\frac{1}{\alpha^2}-\frac{1}{\alpha^2}e^{-\alpha\theta}-\frac{\theta}{\alpha}e^{-\alpha \theta}\right\}+e^{-\alpha\theta}W(0).$
Note that $G(\cdot)$ is continuous on $[0,\infty],$ $G(0)=K+W(0)\ne  W(0)$ because $K>0$, so that $\inf_{\theta\in[0,\infty]}G(\theta)=\inf_{\theta\in(0,\infty]}G(\theta)$, and $G(\infty)=\frac{gh}{\alpha^2}.$ The derivative of $G(\theta)$ in $\theta\in[0,\infty)$ is computed as follows:
\begin{eqnarray*}
&&\frac{dG(\theta)}{d\theta}= -\alpha K e^{-\alpha \theta}+gh\left\{\frac{1}{\alpha}e^{-\alpha \theta}-\frac{1}{\alpha}e^{-\alpha \theta}+\theta e^{-\alpha\theta}\right\}-\alpha e^{-\alpha\theta}W(0)\\
&=&-\alpha K e^{-\alpha \theta}+gh\theta e^{-\alpha \theta}-\alpha e^{-\alpha\theta}W(0)=e^{-\alpha\theta}\left\{-\alpha K+gh\theta-\alpha W(0)\right\}=e^{-\alpha\theta}\tilde{H}(\theta),
\end{eqnarray*}
where $\tilde{H}(\theta)=-\alpha K+gh\theta-\alpha W(0)$. Observe that $\tilde{H}(\theta)$ is strictly increasing in $\theta $ with $\tilde{H}(0)=\lim_{\theta\rightarrow 0}\tilde{H}(\theta)=-\alpha K-\alpha W(0)<0$ because $K>0$ and $W(0)\ge 0$, and $\lim_{\theta\rightarrow \infty}\tilde{H}(\theta)=\infty$. 
Thus, $\inf_{\theta\in[0,\infty]}G(\theta)=\inf_{\theta\in (0,\infty)}G(\theta)$ and $G(\theta)$ admits a unique minimizer $\theta^\ast\in (0,\infty)$ given by the unique stationary point, satisfying 
   \begin{eqnarray}\label{PZ2025Eqn09}
   \theta^\ast=\frac{\alpha(W(0)+K)}{gh};~W(0)=G(\theta^\ast).
   \end{eqnarray}
Therefore, 
\begin{eqnarray*}
&&W(0)=K e^{-\alpha\theta^\ast}+\frac{\alpha W(0)+\alpha K}{\theta^\ast}\left\{\frac{1}{\alpha^2}-\frac{1}{\alpha^2}e^{-\alpha\theta^\ast}-\frac{\theta^\ast}{\alpha}e^{-\alpha\theta^\ast}\right\}+e^{-\alpha \theta^\ast}W(0)\\
&=& K e^{-\alpha \theta^\ast}+\frac{\alpha W(0)}{\theta^\ast}\left\{(1-e^{-\alpha\theta^\ast})\frac{1}{\alpha^2}\right\}+\frac{K}{\theta^\ast}\left\{\frac{1}{\alpha}-\frac{1}{\alpha}e^{-\alpha \theta^\ast}-\theta^\ast e^{-\alpha \theta^\ast}\right\},
\end{eqnarray*}
so that  \begin{eqnarray}\label{PZ2025Eqn08}
W(0)(\alpha \theta^\ast-1+e^{-\alpha\theta^\ast})=\alpha\theta^\ast K e^{-\alpha \theta^\ast}+K\left\{1-e^{-\alpha \theta^\ast}-\alpha\theta^\ast e^{-\alpha \theta^\ast}\right\}.
\end{eqnarray}
Substituting $W(0)=\frac{\theta^\ast gh}{\alpha}-K$ (see (\ref{PZ2025Eqn09})) in the above, we see that
\begin{eqnarray*}
\frac{\theta^\ast gh}{\alpha}(\alpha \theta^\ast-1+e^{-\alpha\theta^\ast})-K(\alpha \theta^\ast-1+e^{-\alpha\theta^\ast})=\alpha\theta^\ast K e^{-\alpha \theta^\ast}+K\left\{1-e^{-\alpha \theta^\ast}-\alpha\theta^\ast e^{-\alpha \theta^\ast}\right\},
\end{eqnarray*}  
i.e., $\frac{\theta^\ast gh}{\alpha}(\alpha \theta^\ast-1+e^{-\alpha\theta^\ast})=K\alpha\theta^\ast$, or, after a rearrangement, 
\begin{eqnarray*}
\frac{gh}{\alpha}\left(1-e^{-\alpha \theta^\ast}\right)=gh\theta^\ast-K\alpha.
\end{eqnarray*}
Standard analysis shows that there is a unique solution in $(0,\infty)$ to the equation \begin{eqnarray*} 
\frac{gh}{\alpha}\left(1-e^{-\alpha x}\right)=ghx-K\alpha.
\end{eqnarray*}
We denote that unique solution in $(0,\infty)$ by $x_g.$ Thus $\theta^\ast=x_g.$ Now from (\ref{PZ2025Eqn08}), we see that
\begin{eqnarray}\label{PZ2025Eqn10}
W(0)=\frac{K(1-e^{-\alpha x_g})}{\alpha x_g-1+e^{-\alpha x_g}}.
\end{eqnarray}

Now consider equation (\ref{eqn6}), i.e., $W(x)=\inf_{\theta\in[0,\infty]}\tilde{G}(x,\theta)$.  Note that $\tilde{G}(x,\theta)$ is continuous in $\theta\in[0,\infty].$
The derivative of function $\tilde{G}(x,\theta)$ is given by 
\begin{eqnarray*}
&&\frac{d\tilde{G}(x,\theta)}{d\theta}=-\alpha K e^{-\alpha \theta}+ghxe^{-\alpha \theta}+gh \left\{\frac{1}{\alpha} e^{-\alpha \theta}-\frac{1}{\alpha}e^{-\alpha \theta}+\theta e^{-\alpha \theta}\right\}-\alpha e^{-\alpha\theta}W(0) \\
&=&-\alpha K e^{-\alpha \theta}+gh x e^{-\alpha \theta}+gh \theta e^{-\alpha\theta}-\alpha e^{-\alpha \theta}W(0)= e^{-\alpha\theta}\left\{-\alpha K +gh x+gh \theta -\alpha W(0)\right\}\\
&=:&  F(x,\theta).
\end{eqnarray*}
At $\theta=0$, $\tilde{G}(x,0)= K+W(0)>0$, and $F(x,0)=-\alpha K+gh x-\alpha W(0)=gh(x-\theta^\ast)=gh(x-x_g)$, by (\ref{PZ2025Eqn09}). 

If $x<x_g$, then $F(x,\theta)<0$ for small $\theta$. Under a fixed $x<x_g$, $F(x,\theta)=0$ at $\theta^\ast(x)=\frac{\alpha(W(0)+K)-ghx}{gh}=\frac{gh\theta^*-ghx}{gh}=\theta^*-x
=x_g-x$ by (\ref{PZ2025Eqn09}),  and $F(x,\theta)>0$ for $\theta>x_g-x.$ Thus, $\theta^\ast(x)=x_g-x$ is the unique minimizer of $\tilde{G}(x,\theta)$, providing the minimum in 
\begin{eqnarray*}
&&W(x)=\inf_{\theta\in[0,\infty]}\tilde{G}(x,\theta)=\tilde{G}(x,\theta^\ast(x))\\
&=& K e^{-\alpha(x_g-x)}+\frac{ghx}{\alpha}\left(1-e^{-\alpha(x_g-x)}\right)+ gh\left\{\frac{1}{\alpha^2}-\frac{1}{\alpha^2}e^{-\alpha (x_g-x)}-\frac{x_g-x}{\alpha}e^{-\alpha(x_g-x)}\right\}\\
&&+e^{-\alpha (x_g-x)}\frac{K(1-e^{-\alpha x_g})}{\alpha x_g-1+e^{-\alpha x_g}},
\end{eqnarray*}
where (\ref{PZ2025Eqn10}) is in use.

If $x= x_g,$ then $F(x,0)=0$, and $F(x,\theta)> 0$ for all $\theta\in (0,\infty)$ because the expression in parentheses in the formula for $F(x,\theta)$ is linear in $\theta$.
Thus,   $\theta^\ast(x)=0=x_g-x$ is the unique minimizer of $\tilde{G}(x,\theta)$, providing the minimum in $\inf_{\theta\in[0,\infty]}\tilde{G}(x,\theta)=W(x)=K+W(0)=K+\frac{K(1-e^{-\alpha x_g})}{\alpha x_g-1+e^{-\alpha x_g}}$, where (\ref{PZ2025Eqn10}) is in use.

If $x>x_g$, then $F(x,\theta)>0$ for all $\theta\in [0,\infty)$, and $\theta^\ast(x)=0$ is the unique minimizer of $\tilde{G}(x,\theta)$, providing the minimum in $\inf_{\theta\in[0,\infty]}\tilde{G}(x,\theta)=W(x)=K+W(0)=K+\frac{K(1-e^{-\alpha x_g})}{\alpha x_g-1+e^{-\alpha x_g}}$.

Note that we have obtained the nonnegative continuous and bounded function $W(\cdot)$ satisfying (\ref{eqn6}). Thus, $W^*_g(\cdot)=W(\cdot)$, Moreover, the minimum in (\ref{eqn6}) is provided by unique minimizer $f(x)=\max\{x_g-x,0\}.$

Second, let $g=0.$ In this case, the nonnegative bounded function $W^*_g(\cdot)=W^*_0(\cdot)$ satisfies equation $W^*_0(x)=\inf_{\theta\in[0,\infty]}\left\{ K e^{-\alpha \theta}+e^{-\alpha \theta}W^*_0(0)\right\}=0$, where the minimum is provided by the unique minimizer $f(x)=\infty.$ 

If we extend $x_g$ by letting  $x_g=\infty$ at $g=0$,
then for all $g\in[0,\infty)$, the unique minimizer providing $W^*_g(x)=\inf_{\theta\in[0,\infty]} \{\bar C_0(x,\theta+g\bar C_1(x,\theta))+\int_{\bf X} W^*_g(y)Q(dy|x,\theta)\}$ is given by the single deterministic stationary strategy $f(x)=\max\{x_g-x,0\}.$
Note that, after the extension, $x_g$, viewed as a function of $g\in[0,\infty)$, provides a bijection from $[0,\infty)$ onto $(0,\infty]$. When $g>0$, $g\in(0,\infty)$ and $x_g\in(0,\infty)$ are in 1-1 correspondence because for any fixed $x>0$ one can put
\begin{equation}\label{e79}
g=\frac{K\alpha^2}{\alpha hx-h+he^{-\alpha x}},
\end{equation}
so that $x=x_g$. Here $g>0$ because the denominator equals zero at $x=0$, and its derivative is $\alpha h-\alpha he^{-\alpha x}>0$ at $x>0$. Clearly, $x_g$ decreases when $g$ increases and $\lim_{g\downarrow 0}x_g=\infty$, $\lim_{g\to\infty}x_g=0$.

We summarize the key observations in this step in the next statement. 
\begin{lemma}\label{t71}
Consider the impulse control problem in Example \ref{PZ2025Example02}. Suppose $g>0$ and let $x_g>0$ be the unique solution to equation
\begin{equation}\label{e78}
\frac{gh}{\alpha}\left(1-e^{-\alpha x}\right)=ghx-K\alpha.
\end{equation}
Then the following assertions hold
\begin{itemize}
\item[(a)] If $x\le x_g$, then
\begin{eqnarray*}
W^*_g(x) &=& Ke^{-\alpha(x_g-x)}+\frac{ghx}{\alpha}\left(1-e^{-\alpha(x_g-x)}\right)\\
&& +gh\left[\frac{1}{\alpha^2}-\frac{1}{\alpha^2}e^{-\alpha(x_g-x)}-\frac{x_g-x}{\alpha}e^{-\alpha(x_g-x)}\right]\\
&& +e^{-\alpha(x_g-x)}\frac{K(1-e^{-\alpha x_g})}{\alpha x_g-1+e^{-\alpha x_g}}.
\end{eqnarray*}
\item[(b)] If $x>x_g$, then
$$W^*_g(x)=K+W^*_g(0)=K+\frac{K(1-e^{-\alpha x_g})}{\alpha x_g-1+e^{-\alpha x_g}}.$$
\end{itemize}
In each case, the single uniformly optimal deterministic stationary strategy $f$ is given by
$$f(x)=\max\{x_g-x,0\}.$$
In words, one has to apply the impulse if and only if the current state $x\ge x_g$.

When $g=0$, $W^*_g(x)\equiv 0$ and $f(x)\equiv\infty$ is the single uniformly optimal deterministic stationary strategy, i.e., no impulses should be applied.

The dual functional $h(\cdot)$ has the form $h(g)=W^*_g(0)-gd=\frac{K(1-e^{-\alpha x_g})}{\alpha x_g-1+e^{-\alpha x_g}}-gd$.
\end{lemma}
  
\textbf{Step 2.} In this step, we find the maximizer of the dual functional  $h(\cdot)$. Since $x_g$ provides a bijection from $[0,\infty)$ onto $(0,\infty]$, as was observed in Step 1, it is convenient to express $h(g)$ in terms of $x_g$ using formulas  (\ref{PZ2025Eqn10}) and (\ref{e79}): one has to maximize function
$$H(x)=\frac{K(1-e^{-\alpha x})}{\alpha x-1+e^{-\alpha x}}-\frac{K\alpha^2 d}{h(\alpha x-1+e^{-\alpha x})},$$
valid also for $x=x_g=\infty$ corresponding to $g=0$.  After that, $g^*$ is calculated using  (\ref{e79}) with $x^*$ being the maximizer of $H(\cdot)$.
  
The derivative of $H(x)$ in $x\in(0,\infty)$ is computed as follows:
\begin{eqnarray*}
\frac{dH(x)}{dx} &=& \frac{h^2K\alpha e^{-\alpha x}(\alpha x-1+e^{-\alpha x})-[hK(1-e^{-\alpha x})-K\alpha^2 d]h(\alpha-\alpha e^{-\alpha x})}{h^2(\alpha x-1+e^{-\alpha x})^2}\\
&=& \frac{hK}{h^2(\alpha x-1+e^{-\alpha x})^2}\{h\alpha^2xe^{-\alpha x}-h\alpha e^{-\alpha x}+h\alpha e^{-2\alpha x}\\
&&-h\alpha+h\alpha e^{-\alpha x}+h\alpha e^{-\alpha x}-h\alpha e^{-2\alpha x}+\alpha^3 d-\alpha^3 de^{-\alpha x}\}\\
&=& \frac{K\alpha}{h(\alpha x-1+e^{-\alpha x})^2}\{\alpha hx e^{-\alpha x}-(1-e^{-\alpha x})(h-d\alpha^2)\}=\frac{K\alpha}{h(\alpha x-1+e^{-\alpha x})^2} D(x),
\end{eqnarray*}
where  $$D(x):=\alpha hx e^{-\alpha x}-(1-e^{-\alpha x})(h-d\alpha^2)=(\alpha h x+(h-\alpha^2 d))e^{-\alpha x}-(h-\alpha ^2 d).$$ 

Observe that $\lim_{x\rightarrow 0}D(x)=0.$
Moreover,
\begin{eqnarray*}
\frac{d D(x)}{dx}&=& \frac{d}{dx} \left\{(\alpha h x+(h-\alpha^2 d))e^{-\alpha x}-(h-\alpha ^2 d)\right\}=\alpha h e^{-\alpha x}-\alpha^2 h x e^{-\alpha x}-\alpha(h-\alpha^2 d)e^{-\alpha x}\\
&=& (h-\alpha h x -(h-\alpha^2 d))e^{-\alpha x}\alpha=\alpha^2 e^{-\alpha x}(\alpha d-hx).
\end{eqnarray*}
Since $\lim_{x\rightarrow 0}\frac{d D(x)}{dx}=\alpha^3d>0$, and   function $\alpha d-hx$ strictly decreases in $x\in(0,\infty)$ with $\lim_{x\rightarrow \infty}(\alpha d-hx)=-\infty$, we see that $D(x)$ increases from $0=\lim_{x\rightarrow 0}D(x)$ for $x>0$ up to its unique stationary point, and, beyond that point, strictly decreases to $\alpha^2d-h=\lim_{x\rightarrow\infty}D(x)$ as $x\rightarrow \infty.$

We now distinguish two cases. 
  
Case (a) $d<\frac{h}{\alpha^2}$. In this case, $\alpha^2d-h=\lim_{x\rightarrow\infty}D(x)<0$, and we see that $\frac{dH(x)}{dx}>0$ for small $x>0$ up to the unique stationary point of $H(\cdot)$, say $x^\ast\in(0,\infty)$, and  $\frac{dH(x)}{dx}<0$ beyond that point. Thus, the unique maximizer of $H(x)$ is given by $x^\ast\in(0,\infty)$, i.e., the unique maximizer of the dual functional $h(g)$ is given by $$g^*=\frac{K\alpha^2}{h(\alpha x^*-1+e^{-\alpha x^*})}\in(0,\infty)$$ in view of (\ref{e79}). Since $\frac{dH(x)}{dx}=0$ if and only if $D(x)=0$, and $H(\cdot)$ has a unique stationary point on $(0,\infty)$, we see that  equation
\begin{eqnarray}\label{e29}
\alpha hx e^{-\alpha x}-(1-e^{-\alpha x})(h-d\alpha^2)=0
\end{eqnarray}
has a unique solution in $(0,\infty)$ given by $x^\ast.$

Case (b)  $d\ge \frac{h}{\alpha^2}$. In this case,  $\alpha^2d-h=\lim_{x\rightarrow\infty}D(x)\ge 0$, and, from the discussion above Case (a), we see that $H(x)$ is nondecreasing, and thus $\infty$ is the maximizer of $H(x)$. In view of the 1-1 correspondence between $g\in[0,\infty)$ and $x_g\in (0,\infty]$, we see that  $g^\ast=0$ is a maximizer of the dual functional $h(g)$.

\textbf{Steps 3 and 4.} In Case (a) of Step 2, i.e., when $d<\frac{h}{\alpha^2}$, we have obtained $g^\ast\in(0,\infty)$. From Step 1, see Lemma \ref{t71}, we see that ${\bf F}$ is a singleton consisting of the deterministic stationary strategy  $f^\ast(x)=\max\{x_{g^*}-x,0\}$, where $x_{g^*}=x^*$ comes from equation (\ref{e29}).  According to Theorem \ref{l8}, $f^\ast$ is necessarily the optimal strategy to the original impulse control problem (\ref{PZZeqn02}). One can also explicitly show that ${\cal V}(0,f^*,\bar{C}_1)=d$:
\begin{eqnarray*}
{\cal V}(0,f^*,\bar C_1)&=& \sum_{i=0}^\infty e^{-\alpha x^* i}\bar C_1(0,x^*)=\frac{h\left[\frac{1}{\alpha^2}-\frac{1}{\alpha^2} e^{-\alpha x^*}-\frac{x^*}{\alpha}e^{-\alpha x^*}\right]}{1- e^{-\alpha x^*}}\\
&=& \frac{h}{\alpha^2}-\frac{hx^*}{\alpha} e^{-\alpha x^*}\cdot\frac{h-d\alpha^2}{\alpha h x^* e^{-\alpha x^*}}=\frac{h}{\alpha^2}-\left[\frac{h}{\alpha^2}-d\right]=d,
\end{eqnarray*}
where the third equality is by (\ref{e29}). Finally,
$${\cal V}(0,f^*,\bar{C}_0)= \sum_{i=0}^\infty e^{-\alpha x^* i}\bar C_0(0,x^*)=\frac{K e^{-\alpha x^*}}{1- e^{-\alpha x^*}}.$$

In Case (b) of Step 2, i.e., when $d\ge \frac{h}{\alpha^2}$, we have obtained $g^\ast=0$. From Step 1, see Lemma \ref{t71}, ${\bf F}$ is again a singleton consisting of the deterministic stationary strategy $f^\ast(x)\equiv \infty$. Similar argument as above shows that $f^*$ is the solution to original impulse control problem (\ref{PZZeqn02}) and
$${\cal V}(0,f^*,\bar{C}_0)=\bar C_0(0,\infty)=0,~~~{\cal V}(0,f^*,\bar{C}_1)=\bar C_1(0,\infty)=\frac{h}{\alpha^2}.$$

We now summarize the last conclusions in the following statement. 
\begin{theorem}\label{t72}
Consider the impulse control problem in Example \ref{PZ2025Example02}.\begin{itemize}  
\item[(a)] If $d<\frac{h}{\alpha^2}$, then
$$g^*=\frac{K\alpha^2}{h(\alpha x^*-1+e^{-\alpha x^*})}$$ maximizes the dual  functional $h(g)$,
where $x^*>0$ is the unique positive solution to equation
\begin{eqnarray*} 
(1-e^{-\alpha x})(h-d\alpha^2)=\alpha h x e^{-\alpha x}.
\end{eqnarray*}
The optimal solution to problem (\ref{PZZeqn02}) for the impulse control problem in Example \ref{PZ2025Example02} is given by the deterministic stationary strategy
$$f^*(x)=\max\{x^*-x,0\}:$$
one has to apply the impulse if and only if the current state $x\ge x^*$, with
$${\cal V}(0,f^*,\bar{C}_0)=\frac{Ke^{-\alpha x^*}}{1- e^{-\alpha x^*}};~~~{\cal V}(0,f^*,\bar{C}_1)=d.$$
\item[(b)] If $d\ge\frac{h}{\alpha^2}$, then $g^*=0$ maximizes the dual functional $h(g)$, and the optimal solution to problem (\ref{PZZeqn02}) for the impulse control problem in Example \ref{PZ2025Example02}  is given by the deterministic stationary strategy $f^*(x)\equiv \infty$: no impulses should be applied, with
$${\cal V}(0,f^*,\bar{C}_0)=0;~~~{\cal V}(0,f^*,\bar{C}_1)=\frac{h}{\alpha^2}.$$
\end{itemize}
\end{theorem}

\section{Proofs of the statements}\label{secap}
\par\noindent\textit{Proof of Proposition \ref{pr1}.} $\bf B$ is compact (in the product topology) as the product of compact sets. Thus, (a) in Definition \ref{PZ2025Def01} is verified.

If $u:~{\bf X}_\Delta\to\RR$ is a bounded continuous function, then
$$\int_{{\bf X}_\Delta} u(y)Q(dy|x,b=(\theta,a))=\left\{\begin{array}{ll}
e^{-\alpha\theta}u(l(\phi(x,\theta),a))\\
+(1-e^{-\alpha\theta})u(\Delta) & \mbox{ if } x\ne\Delta,~\theta\ne+\infty;\\
u(\Delta) & \mbox{ otherwise}\end{array}\right.$$
is a continuous function in $(x,b)\in{\bf X}_\Delta\times{\bf B}$. Thus, (b) in Definition \ref{PZ2025Def01} is verified.

Let us verify (c) Definition \ref{PZ2025Def01}, that is, we show for a fixed $j\in\{0,\dots,J\}$ that $\bar C_j(\cdot,\cdot)$ is lower semicontinuous.
Its lower semicontinuity on $\{\Delta\}\times {\bf B}$ follows from $\bar{C}_j(\Delta,b)\equiv 0$ and that $\Delta$ is isolated. Now it remains to show that $\bar{C}_j(\cdot,\cdot)$ is lower semicontinuous on ${\bf X}\times {\bf B}$, as follows.
Recall that a non-negative function on ${\bf X}\times {\bf B}$ is lower semicontinuous if and only if it is the (pointwise) limit of an increasing sequence of bounded non-negative continuous functions, see \cite[Lemma 7.14]{Bertsekas:1978}. Thus, using the monotone convergence theorem, we see that the function
\begin{eqnarray*}
\int_0^\theta e^{-\alpha t}  C^g_j(\phi(x,t))dt= \lim_{m\to\infty}  \int_0^\theta e^{-\alpha t}  C^{g,m}_j(\phi(x,t))dt
\end{eqnarray*}
is non-negative and lower semicontinuous in $(x,b=(\theta,a))\in {\bf X}\times{\bf B}$. Here $C^{g,m}_j(\cdot)\uparrow C^g_j(\cdot)$  pointwise as $m\to\infty$ for non-negative, bounded and continuous functions $C^{g,m}_j(\cdot)$ on ${\bf X}$, whose existence follows from the lower semicontinuity and non-negativity of the function $C^g_j(\cdot)$ and the aforementioned fact. On the other hand, the function defined by $e^{-\alpha\theta} C^I_j(\phi(x,\theta),a)$ is lower semicontinuous on ${\bf X}\times[0,\infty)\times{\bf A}$, and
$$\liminf_{n\to\infty} e^{-\alpha\theta_n}C^I_j(\phi(x_n,\theta_n),a_n)\ge 0=e^{-\alpha\infty}C^I_j(\phi(x,\theta),a)$$
if $(x_n,(\theta_n,a_n))\to (x,(\infty,a))\in{\bf X}\times \bf{B}$.
\hfill $\Box$
\bigskip

\par\noindent\textit{Proof of Lemma \ref{l4}.} Consider a strategy $\pi$ such that $\mu^\pi({\bf X}\times{\bf B})=\infty.$ Recall that $\mu^\pi$ satisfies 
\begin{eqnarray*}
\mu^\pi(dx\times {\bf B})=\delta_{x_0}(dx)+\int_{\bf X}\int_{[0,\infty)}\int_{\bf A}e^{-\alpha\theta}\delta_{l(\phi((y,\theta),a))}(dx)\mu^{\pi}(dy\times d\theta\times da)
\end{eqnarray*}
on $({\bf X},{\cal B}({\bf X})).$ Since $\mu^\pi({\bf X}\times{\bf B})=\infty$, we have $$\infty=1+\int_{\bf X}\int_{[0,\infty)}\int_{\bf A}e^{-\alpha\theta} \mu^{\pi}(dy\times d\theta\times da)=1+ \int_{[0,\infty)}e^{-\alpha\theta} \mu^{\pi}({\bf X}\times d\theta\times {\bf A}).$$ Now, 
\begin{eqnarray*}
{\cal V}(x_0,\pi,\bar{C}_0)&=& \int_{{\bf X}_\Delta\times [0,\infty]\times{\bf A}} \bar{C}_0(x,(\theta,a))\mu^\pi(dx\times d\theta\times da)\\
&\ge& \int_{[0,\infty)} e^{-\alpha\theta} \delta\mu^\pi({\bf X}\times d\theta\times {\bf A})=\infty,
\end{eqnarray*}
as required. $\hfill\Box$
\bigskip

\par\noindent\textit{Proof of Lemma \ref{l7}.} (a) Consider the deterministic stationary strategy identified by $f(x)\equiv (\infty,\hat a)$ with an arbitrarily fixed $\hat a\in{\bf A}$. Under Condition \ref{con5}(b), ${\cal V}(x,f, \bar{C}_0+\sum_{j=1}^J g_j \bar{C}_j)=\EE^f_x\left[\sum_{i=1}^\infty\left(\bar C_0(X_{i-1},B_i)+\sum_{j=1}^J g_j\bar C_j(X_{i-1},B_i)\right)\right]  \le(J+1){\cal C}\int_0^\infty e^{-\alpha t}dt$ for all $x\in {\bf X}_\Delta$. Thus, the function ${\cal V}^\ast(\cdot, \bar{C}_0+\sum_{j=1}^J g_j \bar{C}_j)$  is non-negative and bounded by $(J+1){\cal C}\int_0^\infty e^{-\alpha t}dt$. Since the induced MDP is semicontinuous by assumption, according to  Corollary 9.17.2 of \cite{Bertsekas:1978}, we see that ${\cal V}^\ast(\cdot, \bar{C}_0+\sum_{j=1}^J g_j \bar{C}_j)$ is a bounded, non-negative and lower semicontinuous function on ${\bf X}_\Delta$. It follows that $W^\ast_{\bar g}(\cdot)$, defined as the restriction of ${\cal V}^\ast(\cdot, \bar{C}_0+\sum_{j=1}^J g_j \bar{C}_j)$ on ${\bf V}$, is bounded, nonnegative and lower semicontinuous on ${\bf V}$.  It is clear that ${\cal V}^\ast(x, \bar{C}_0+\sum_{j=1}^J g_j \bar{C}_j)(x)=0$ for all $x\in {\bf V}^c$. It follows from this fact and Proposition 9.8 of \cite{Bertsekas:1978} that 
\begin{eqnarray}\label{PZ2025Eqn06}
&&{\cal V}^\ast(x, \bar{C}_0+\sum_{j=1}^J g_j \bar{C}_j)\nonumber\\
 &=& \inf_{b\in{\bf B}}\left\{\bar C_0(x,b)+\sum_{j=1}^J g_j\bar C_j(x,b)+\int_{{\bf X}_{\Delta}} {\cal V}^\ast(y, \bar{C}_0+\sum_{j=1}^J g_j \bar{C}_j) Q(dy|x,b)\right\}\nonumber\\
&=&\inf_{b\in{\bf B}}\left\{\bar C_0(x,b)+\sum_{j=1}^J g_j\bar C_j(x,b)+\int_{\bf V} {\cal V}^\ast(y, \bar{C}_0+\sum_{j=1}^J g_j \bar{C}_j) Q(dy|x,b)\right\}
\end{eqnarray}
for all $x\in{\bf X}_{\Delta}.$ Hence,  $W^\ast_{\bar g}(\cdot)$, defined as the restriction of ${\cal V}^\ast(\cdot, \bar{C}_0+\sum_{j=1}^J g_j \bar{C}_j)$ on ${\bf V}$, satisfies (\ref{e28}), and is feasible in problem   (\ref{e27})-(\ref{e28}).

 (b) According to Proposition 9.12 of \cite{Bertsekas:1978}, a deterministic stationary strategy $\hat{f}$ is uniformly optimal if and only if it attains the infimum in  equality (\ref{PZ2025Eqn06}) for all $x\in {\bf X}_{\Delta}$. For $x\in {\bf V}^c$, the infimum is evidently attained at $f^\ast(x)$. Consequently, (b) is proved once we show that $W(\cdot)=W^\ast_{\bar{g}}(\cdot)$ for the given function $W(\cdot)\in\mathbb{B}({\bf V})$ satisfying equalities (\ref{enn})
for a measurable mapping $\hat{f}(\cdot):~{\bf V}\to{\bf B}$. Note that if $W(\cdot)\in\mathbb{B}({\bf V})$ is lower semicontinuous, then such a measurable mapping $\hat{f}(\cdot):~{\bf V}\to{\bf B}$ exists, see e.g., Theorem 3.3 of \cite{Feinberg2013JMAA}.   It is convenient to extend $W(\cdot)$ to ${\bf X}$ by putting $W(x)=0$ for all $x\in {\bf V}^c$ in the rest of this proof.

 First, we show that  $W(x)\le{\cal V}(x,\pi,\bar C_0)+\sum_{j=1}^J g_j{\cal V}(x,\pi,\bar C_j)$ for all $x\in {\bf V}$ and strategies $\pi$. Let $x\in {\bf V}$ and a strategy $\pi$ be fixed.
 We distinguish two cases. 
 
 Case 1. ${\cal V}(x,\pi,\bar{C}_0)=\infty$. In this case, the desired inequality holds automatically.
 
 Case 2. ${\cal V}(x,\pi,\bar{C}_0)<\infty$. Let us show in this case that 
$T_\infty:=\lim_{i\to\infty}T_i=\infty$ $\PP^\pi_x$-a.s..  Suppose for contradiction that  $\PP^\pi_x(T_\infty<\infty)>0.$ Since, $\{T_\infty<\infty\}=\bigcup_{n=1}^\infty\{T_\infty<n\}$, there is some $T\in\{1,2,\dots\}$ such that 
$\PP_x^\pi(T_\infty<T)>0.$ Since $\{T_\infty<T\}\subseteq \bigcap_{i=1}^\infty\{T_i<T\},$ it holds that $\lim_{i\rightarrow \infty}\PP_x^\pi(T_i<T)>0$. Now
\begin{eqnarray*}
&&{\cal V}(x,\pi,\bar{C}_0)\ge \sum_{i=1}^\infty \EE_x^\pi\left[e^{-\alpha T_i}\delta\right]\ge  \sum_{i=1}^\infty \EE_x^\pi\left[ \II\{T_i<T\}e^{-\alpha T_i}\delta\right]\\
&\ge& \delta \sum_{i=1}^\infty \EE_x^\pi\left[\II\{T_i<T\}e^{-\alpha T}\right]=\delta e^{-\alpha T} \sum_{i=1}^\infty \PP_x^\pi(T_i<T)=\infty,
\end{eqnarray*}
where the last equality holds because, for each $i=1,2,\ldots$, $\PP_x^\pi(T_i<T)\ge\lim_{i\rightarrow \infty}\PP_x^\pi(T_i<T)>0$. This yields a contradiction, proving that $T_\infty=\lim_{i\to\infty}T_i=\infty$ $\PP^\pi_x$-a.s.. Now,  since $\PP^\pi_x(X_I\in{\bf X})=\EE^\pi_x[e^{-\alpha T_I}]$,
\begin{eqnarray*}
&&\EE^\pi_x[|W(X_I)|]=\EE_x^\pi \left[\II\{X_I\in {\bf V}\} |W(X_I)|\right]\le \EE_x^\pi \left[\II\{X_I\in {\bf V}\}\right] \sup_{y\in {\bf V}}|W(y)|\\
&\le& \PP_x^\pi (X_I\in {\bf X}) \sup_{y\in {\bf V}}|W(y)|=\EE_x^\pi\left[e^{-\alpha T_I}\right]\sup_{y\in {\bf V}}|W(y)|\to 0 \mbox{ as } I\to\infty,
\end{eqnarray*}
where the convergence to $0$ holds because of the dominated convergence theorem and the fact that $T_\infty=\lim_{i\to\infty}T_i=\infty$ $\PP^\pi_x$-a.s.. 
Furthermore, for each $I=1,2,\ldots$,
\begin{eqnarray}\label{PZ2025Eqn07}
0 &\le & \EE^\pi_x\left[\sum_{i=1}^I\left\{\bar C_0(X_{i-1},B_i)+\sum_{j=1}^J g_j\bar C_j(X_{i-1},B_i)-W(X_{i-1})\right.\right.\nonumber\\
&&\left.\left.\vphantom{\sum_i^I}+\int_{\bf V} W(y)Q(dy|X_{i-1},B_i)\right\}\right]\nonumber\\
&=&\EE^\pi_x\left[\sum_{i=1}^I\left\{\bar C_0(X_{i-1},B_i)+\sum_{j=1}^J g_j\bar C_j(X_{i-1},B_i)\right\}\right]-\EE^\pi_x[W(X_0)]\nonumber\\
&&+\EE^\pi_x[W(X_1)]-\EE^\pi_x[W(X_1)]+\EE^\pi_x[W(X_2)]-\ldots+\EE^\pi_x[W(X_I)],
\end{eqnarray}
where for the last equality we recall that $W(x)=0$ for  $x\in {\bf V}^c.$
Note that all the terms here are finite, because the functions $\bar{C}_j(\cdot,\cdot)$ and $W(\cdot)$ are all bounded.
After passing to the limit as $I\to\infty$, we obtain the claimed inequality
$$W(x)\le {\cal V}(x,\pi,\bar C_0)+\sum_{j=1}^J g_j{\cal V}(x,\pi,\bar C_j),$$ too.

Secondly, we show that $W(x)= {\cal V}(x,\hat{f},\bar C_0)+\sum_{j=1}^J g_j{\cal V}(x,\hat{f},\bar C_j)$ for each $x\in {\bf V}$, as follows.

Let $x\in {\bf V}$ be fixed again. For the strategy $\pi=\hat{f}$ in the statement of this lemma,  we have equalities in (\ref{PZ2025Eqn07}):
\begin{eqnarray*}0&=&\EE^{\hat{f}}_x\left[\sum_{i=1}^I\left\{\bar C_0(X_{i-1},B_i)+\sum_{j=1}^J g_j\bar C_j(X_{i-1},B_i)\right\}\right]-\EE^\pi_x[W(X_0)]+\EE^{\hat{f}}_x[W(X_I)].
\end{eqnarray*}
Since $W(x)$ and $\limsup_{I\to\infty} \EE^{\hat{f}}_x[W(X_I)]$ are finite, upon passing to the limit as $I\rightarrow \infty,$ we see that $0= {\cal V}(x,\hat{f},\bar{C}_0)+\sum_{j=1}^J g_j{\cal V}(x,\hat{f},\bar{C}_j)-W(x)+\limsup_{I\to\infty} \EE^{\hat{f}}_x[W(X_I)],$
i.e., $W(x)={\cal V}(x,\hat{f},\bar{C}_0)+\sum_{j=1}^J g_j{\cal V}(x,\hat{f},\bar{C}_j)+\limsup_{I\to\infty} \EE^{\hat{f}}_x[W(X_I)].$ Again, since $W(x)$ and $\limsup_{I\to\infty} \EE^{\hat{f}}_x[W(X_I)]$ are finite, we see that ${\cal V}(x,\hat{f},\bar{C}_0)+\sum_{j=1}^J g_j{\cal V}(x,\hat{f},\bar{C}_j)<\infty$. In particular,
${\cal V}(x,\hat{f},\bar{C}_0)<\infty$, and we are in Case 2 above with $\pi=\hat{f}$. Applying the result established in Case 2 for $\pi=\hat{f}$, we see that $\lim_{I\to\infty} \EE^{\hat{f}}_x[W(X_I)]=0$. Hence,
$$W(x)= {\cal V}(x,\hat{f},\bar{C}_0)+\sum_{j=1}^J g_j{\cal V}(x,\hat{f},\bar{C}_j).$$

Finally, combining the results established in the above two steps, we see that the statement of the lemma holds. \hfill$\Box$
\bigskip
 
\par\noindent\textit{Proof of Theorem \ref{t5}.} (a) Let $\bar g\in\RR_+^J$ be arbitrarily fixed.  In view of Lemma \ref{l7}, we only need to show that the bounded, non-negative and lower semicontinuous function $W^\ast_{\bar{g}}(\cdot)$ defined by (\ref{e31})
solves problem (\ref{e27})-(\ref{e28}) with $val(\mbox{problem (\ref{e27})-(\ref{e28})})=W^*_{\bar g}(x_0)$. This is done as follows, by making use of the equivalence (up to an additive constant) between problem (\ref{e27})-(\ref{e28}) and problem (\ref{e26}). 

Since the induced MDP is semicontinuous by assumption, according to  Corollary 9.17.2 of \cite{Bertsekas:1978},  there exists a uniformly optimal deterministic stationary strategy $\hat f$ satisfying
${\cal V}^\ast(x, \bar{C}_0+\sum_{j=1}^J g_j \bar{C}_j)= \inf_\pi {\cal V}(x,\pi, \bar{C}_0+\sum_{j=1}^J g_j \bar{C}_j)=  {\cal V}(x,\hat{f}, \bar{C}_0+\sum_{j=1}^J g_j \bar{C}_j)$ for all $x\in {\bf X}_\Delta$. Without loss of generality, we put $\hat f(x)=f^*(x)$ for $x\in{\bf V}^c$. (Recall equalities (\ref{e5}).)  

Consider the following problem, which is dual with problem  (\ref{e26}):
 \begin{equation}\label{e32}
\mbox{Minimize over $\mu\in{\cal M}^{<\infty}({\bf V}\times {\bf B})$: }\sup_{W(\cdot)\in \mathbb{B}({\bf V})} L_2(\mu,W,\bar g).
\end{equation}
We know that
\begin{eqnarray}\label{PZ2025Eqn05}
val(\mbox{problem (\ref{e32})})\ge val(\mbox{problem (\ref{e26})}).
\end{eqnarray}
(See, e.g., (1.8) in \cite{rock}.)
 
On the other hand, problem (\ref{e32}), up to the constant $\sum_{j=1}^J g_jd_j$, is equivalent to the following one:
\begin{eqnarray*}\mbox{Minimize over }\mu\in{\cal M}^{<\infty}({\bf V}\times{\bf B})&:&\int_{{\bf V}\times{\bf B}} \left(\bar C_0(x,b)+\sum_{j=1}^J g_j\bar C_j(x,b)\right)\mu(dx\times db)\\
\mbox{ subject to}&:&\mu(\Gamma\times \textbf{B})=\delta_{x_0}(\Gamma)+\int_{{\bf V}\times \textbf{B}}Q(\Gamma|y,b)\mu(dy\times db),\\
&&\Gamma\in{\cal B}({\bf V}),
\end{eqnarray*}
where we recall that the above constraint equality is the same as  (\ref{e22}).
Indeed,
if for some $\mu\in {\cal M}^{<\infty}({\bf V}\times{\bf B})$ equality  (\ref{e22}) is not satisfied, then, similarly to Remark \ref{r2},  $\sup_{W(\cdot)\in \mathbb{B}({\bf V})} L_2(\mu,W,\bar g)=+\infty$; otherwise, $$L_2(\mu,W,\bar g)+\sum_{j=1}^Jg_jd_j=\int_{{\bf V}\times{\bf B}} \left(\bar C_0(x,b)+\sum_{j=1}^J g_j\bar C_j(x,b)\right)\mu(dx\times db)$$  for all $W(\cdot)\in \mathbb{B}({\bf V})$.

Now consider the occupation measure $\mu^{\hat{f}}$ of a uniformly optimal deterministic stationary strategy $\hat{f}$. Since  ${\cal V}(x_0,\hat{f}, \bar{C}_0)\le {\cal V}(x_0,\hat{f}, \bar{C}_0+\sum_{j=1}^J g_j \bar{C}_j)\le (J+1){\cal C}\int_0^\infty e^{-\alpha t}dt$, where the last inequality follows from the uniform optimality of $\hat{f}$ and the rightmost bound was observed in the proof of Lemma \ref{l7}(a), we see from Lemma \ref{l4} that $\mu^{\hat{f}}({\bf V}\times {\bf B})\le \mu^{\hat{f}}({\bf X}\times {\bf B})<\infty.$  Consider $\mu_{\bar g}:=\mu^{\hat f}|_{{\bf V}\times {\bf B}}$. Then $\mu_{\bar g}\in {\cal M}^{<\infty}({\bf V}\times {\bf B})$ and satisfies equality (\ref{e22}) because $\hat f(x)=f^*(x)$ for $x\in{\bf V}^c$. (Recall (\ref{e5}.) Moreover, it holds that 
$$\int_{{\bf V}\times{\bf B}} \left(\bar C_0(x,b)+\sum_{j=1}^J g_j\bar C_j(x,b)\right)\mu_{\bar g}(dx\times db)=W^*_{\bar g}(x_0).$$ Therefore, 
\begin{eqnarray}\label{enum22}
val(\mbox{problem (\ref{e32})})&\le& \int_{{\bf V}\times{\bf B}} \left(\bar C_0(x,b)+\sum_{j=1}^J g_j\bar C_j(x,b)\right)\mu_{\bar g}(dx\times db)-\sum_{j=1}^Jg_jd_j\nonumber\\
&=&W^*_{\bar g}(x_0)-\sum_{j=1}^Jg_jd_j.
\end{eqnarray}

Since problem  (\ref{e26}) is,  up to the additive constant $\sum_{j=1}^Jg_jd_j$, equivalent to (\ref{e27})-(\ref{e28}), and $W^*_{\bar g}(\cdot)$ is a feasible solution in problem  (\ref{e27})-(\ref{e28}), as was observed in Lemma \ref{l7}(a), we see that
$$val(\mbox{problem (\ref{e26})})\ge W^*_{\bar g}(x_0)-\sum_{j=1}^Jg_jd_j.$$
Thus, by (\ref{enum22}),
$$val(\mbox{problem (\ref{e26})})\ge val(\mbox{problem (\ref{e32})}), $$
which, together with (\ref{PZ2025Eqn05}), leads to 
$$val(\mbox{problem (\ref{e26})})= val(\mbox{problem (\ref{e32})})=W^*_{\bar g}(x_0)-\sum_{j=1}^Jg_jd_j.$$
Hence, $W^*_{\bar g}(\cdot)$ is an optimal solution in problem (\ref{e26}). It follows that  $W^*_{\bar g}(\cdot)$ solves problem (\ref{e27})-(\ref{e28}), which is equivalent to problem (\ref{e26}) up to an additive constant.

(b) Consider the strategy $\hat{f}$ as in the proof of part (a), and adopt the notations introduced therein. Then 
\begin{eqnarray*}
&&W^{\ast}_{\bar{g}}(x_0)-\sum_{j=1}^Jg_jd_j=\inf_{\pi} \EE^\pi_{x_0}\left[\sum_{i=1}^\infty\left( \bar C_0(X_{i-1},B_i)+\sum_{j=1}^Jg_j\bar C_j(X_{i-1},B_i)\right)\right]-\sum_{j=1}^J g_jd_j\\
&=&{\cal V}(x_0,\hat{f},\bar{C}_0+\sum_{j=1}^Jg_j \bar{C}_j)-\sum_{j=1}^J g_jd_j\\
&=&\int_{{\bf V}\times {\bf B}}\left(\bar{C}_0(x,b)+\sum_{j=1}^J g_j\bar{C}_j(x,b)\right)\mu^{\hat{f}}(dx\times db)-\sum_{j=1}^J g_jd_j.
\end{eqnarray*} 
We have seen in the proof of part (a) of this theorem that $\mu^{\hat{f}}({\bf X}\times {\bf B})<\infty.$  The non-negative functions $\bar{C}_j(\cdot,\cdot)$, $j\in\{0,\dots,J\}$, are bounded from above by $(\frac{1}{\alpha}+1){\cal C}<\infty.$ It follows that 
\begin{eqnarray*}
\int_{{\bf V}\times{\bf B}}\bar{C}_j(x,b)\mu^{\hat{f}}(dx\times db)\le  (\frac{1}{\alpha}+1){\cal C}\mu^{\hat{f}}({\bf V}\times {\bf B})<\infty
\end{eqnarray*}
for each $j\in\{0,\dots,J\}$, and $\mu^{\hat{f}}|_{{\bf V}\times {\bf B}} \in {\cal D}$. One may write 
$$\int_{{\bf V}\times {\bf B}}\left(\bar{C}_0(x,b)+\sum_{j=1}^J g_j\bar{C}_j(x,b)\right)\mu^{\hat{f}}(dx\times db)-\sum_{j=1}^J g_jd_j=L_1(\mu^{\hat{f}}|_{{\bf V}\times{\bf B}},\bar{g}).$$ Hence,
\begin{eqnarray*}
\inf_{\pi} \EE^\pi_{x_0}\left[\sum_{i=1}^\infty\left( \bar C_0(X_{i-1},B_i)+\sum_{j=1}^Jg_j\bar C_j(X_{i-1},B_i)\right)\right]-\sum_{j=1}^J g_jd_j\ge\inf_{\mu\in{\cal D}} L_1(\mu,\bar g)= h(\bar g).
\end{eqnarray*}
Next, we prove that the opposite direction of the above inequality holds, too.
For each $\mu\in {\cal D}$, consider its induced strategy $\pi^\mu$. By  Lemma \ref{l3}, we see that  
\begin{eqnarray*}
&&\int_{{\bf V}\times {\bf B}}\left(\bar{C}_0(x,b)+\sum_{j=1}^J g_j\bar{C}_j(x,b)\right)\mu(dx\times db)\\
&\ge& \EE^{\pi^\mu}_{x_0}\left[\sum_{i=1}^\infty\left( \bar C_0(X_{i-1},B_i)+\sum_{j=1}^Jg_j\bar C_j(X_{i-1},B_i)\right)\right]
\end{eqnarray*} 
and thus,
$L_1(\mu,\bar g)\ge \EE^{\pi^\mu}_{x_0}\left[\sum_{i=1}^\infty\left( \bar C_0(X_{i-1},B_i)+\sum_{j=1}^Jg_j\bar C_j(X_{i-1},B_i)\right)\right]-\sum_{j=1}^J g_jd_j.$ Since $\mu\in {\cal D}$ was arbitrarily fixed, it follows that
  $$h(\bar g)=\inf_{\mu\in{\cal D}} L_1(\mu,\bar g)\ge \inf_{\pi} \EE^\pi_{x_0}\left[\sum_{i=1}^\infty\left( \bar C_0(X_{i-1},B_i)+\sum_{j=1}^Jg_j\bar C_j(X_{i-1},B_i)\right)\right]-\sum_{j=1}^J g_jd_j.$$ The first assertion in part (b) is thus proved, whereas the second assertion follows from the first assertion, part (a) of this theorem and Theorem \ref{prop2}(a).

(c) This assertion follows from parts (a) and (b).

(d) The assertion in this part follows from the following chain of equalities. According to part (b),
$$\sup_{(W(\cdot),\bar g)\in \mathbb{B}({\bf V})\times\RR_+^J} \inf_{\mu\in {\cal M}^{<\infty}({\bf V}\times {\bf B})} L_2(\mu,W,\bar g)=val(\mbox{problem (\ref{e25})})=h(\bar g^*).$$
By Theorem \ref{prop2}(a),
$h(\bar g^*)=val(\mbox{problem (\ref{enum4})})=val(\mbox{primal convex program (\ref{enum3})}).$
As was noted at the beginning of Subsection \ref{PZ2025Subsec01},
$$val(\mbox{primal convex program (\ref{enum3})})=val(\mbox{relaxed problem (\ref{SashaLp02})-(\ref{e10})}).$$
Finally, under the imposed conditions,
\begin{eqnarray*}
&&val(\mbox{relaxed problem (\ref{SashaLp02})-(\ref{e10})})=val(\mbox{second convex program (\ref{e21})-(\ref{e23})})\\
&=&val(\mbox{problem (\ref{e24})})= \inf_{\mu\in{\cal M}^{<\infty}({\bf V}\times {\bf B})}\sup_{(W(\cdot),\bar g)\in \mathbb{B}(\bf V)\times\RR_+^J} L_2(\mu,W,\bar g),
\end{eqnarray*}
where the first equality was observed in the paragraph below Lemma \ref{l4}, and the second equality holds by Remark \ref{r2}.
\hfill$\Box$
\bigskip

\par\noindent\textit{Proof of Lemma \ref{l6}.} By Corollary 9.17.2 of \cite{Bertsekas:1978}, the Bellman function ${\cal V}^\ast(\cdot,C)$ is lower semicontinuous and there exists a uniformly optimal deterministic stationary strategy $\breve f$ in the sense of (\ref{e11}).  Without loss of generality we accept that $\breve f(\Delta)=f^*(\Delta)$.  Let us show that 
\begin{equation}\label{e12}
{\cal V}(y,f,C):=
\EE^f_{y}\left[\sum_{i=1}^\infty C(Y_{i-1},f(Y_{i-1}))\right]={\cal V}^\ast(y,C)
\end{equation}
for $\mu^f(dy\times{\bf B})$-almost all $y\in{\bf V}$.

Suppose for contradiction that there is $\Gamma_Y\in{\cal B}({\bf V})$ such that equality (\ref{e12}) is violated for all $y\in\Gamma_Y$ and $\mu^f(\Gamma_Y\times{\bf B})>0$. Let $I$ be the first moment when $\EE^f_{y_0}[\II\{Y_I\in\Gamma_X\}]>0$. Then
\begin{eqnarray*}
\infty>\EE^f_{y_0}\left[\sum_{i=1}^\infty C(Y_{i-1},f(Y_{i-1}))\right]&=&\EE^f_{y_0}\left[\sum_{i=1}^I C(Y_{i-1},f(Y_{i-1}))\right]\\
&&+\EE^f_{y_0}\left[\EE^f_{Y_I}\left[\sum_{i=1}^\infty C(Y_{i-1},f(Y_{i-1}))\right]\right],
\end{eqnarray*}
where the first inequality holds because ${\cal V}^\ast(y_0,C)<\infty.$
The last term in the above equality is strictly larger than $\EE^f_{y_0}[{\cal V}^*(Y_I,C)]$. At the same time, for the  deterministic Markov strategy 
$$\hat f_i(y):=\left\{\begin{array}{ll} f(y), & \mbox{ if } i\le I,\\ \breve f(y), & \mbox{ if } i>I\end{array}\right.$$ 
we have
\begin{eqnarray*}\EE^{\hat f}_{y_0}\left[\sum_{i=1}^\infty C(Y_{i-1},f(Y_{i-1}))\right]&=&\EE^f_{y_0}\left[\sum_{i=1}^I C(Y_{i-1},f(Y_{i-1}))\right]
+\EE^f_{y_0}\left[{\cal V}^*(Y_I,C)\right]\\
&<&\EE^f_{y_0}\left[\sum_{i=1}^\infty C(Y_{i-1},f(Y_{i-1}))\right],
\end{eqnarray*}
and thus the strategy $f$ is not optimal at $x_0$, yielding a contradiction. This proves equality (\ref{e12}).

Now we modify the strategy $f$ on the $\mu^f(dy\times{\bf B})$-null subset of $\bf V$:
$$\tilde f(y):=\left\{\begin{array}{ll} f(y), & \mbox{ if (\ref{e12}) holds for $y$},\\ \breve f(y), & \mbox{ otherwise.} \end{array}\right.$$
The left-hand side in (\ref{e12}) is a measurable function of $y$: see Propositions C.1.1(b), B.1.31 and B.1.34 of \cite{book20}. Thus, $\tilde f(\cdot)$ is a measurable mapping because the set
$\{y\in{\bf V}:~(\ref{e12})\mbox{ holds for }y\}$
is measurable.

Next, let us prove the following equality
\begin{equation}\label{enum0}
{\cal V}^\ast(y,C)=C(y,\tilde f(y))+\int_{\bf V} {\cal V}^\ast(z,C)Q(dz|y,\tilde f(y)),~~~~~y\in{\bf V}.
\end{equation}
Let $y\in {\bf V}$ be fixed. If  equality (\ref{e12}) is not valid at this $y\in {\bf V}$, then $\tilde f(y)=\breve f(y)$, and equality (\ref{enum0}) holds by Proposition 9.12 of \cite{Bertsekas:1978}. Suppose equality (\ref{e12}) is valid for this $y\in{\bf V}$, so that $\tilde f(y)=f(y)$. Then
$${\cal V}^\ast(y,C)={\cal V}(y,f,C)=C(y,f(y))+\int_{\bf V}{\cal V}(z,f,C)Q(dz|y,f(y)),$$
where
$${\cal V}(z,f,C)=\EE^f_z\left[\sum_{i=1}^\infty C(Y_{i-1},B_i)\right]\ge {\cal V}^\ast(z,C).$$
Thus,
$${\cal V}^\ast(y,C)\ge C(y,f(y))+\int_{\bf V} {\cal V}^\ast(z,C)Q(dz|y, f(y)).$$
If this inequality was strict, then we would have
$${\cal V}^\ast(y,C)>\inf_{b\in{\bf B}}\left\{C(y,b)+\int_{\bf V}{\cal V}^\ast(z,C)Q(dz|y,b)\right\},$$
which is impossible, because the Bellman function ${\cal V}^\ast(\cdot,C)$ satisfies the Bellman equation
$$W(y)=\inf_{b\in{\bf B}}\left\{C(y,b)+\int_{\bf V} W(z)Q(dz|y,b)\right\},~y\in{\bf V}.$$ Consequently,
we conclude that equality (\ref{enum0}) is valid for all $y\in{\bf V}$, and the deterministic stationary strategy $\tilde f$ is uniformly optimal, again, by Proposition 9.12 of \cite{Bertsekas:1978}.

To show that the occupation measure $\mu^{\tilde f}$ coincides with $\mu^f$ on ${\bf V}\times{\bf B}$, firstly note that
\begin{equation}\label{e55}
\mu^f(dy\times db)=\mu^f(dy\times{\bf B})\delta_{f(y)}(db)=\mu^f(dy\times{\bf B})\delta_{\tilde f(y)}(db)
\end{equation}
because $f(\cdot)=\tilde f(\cdot)$ for $\mu^f(dy\times{\bf B})$-almost all $y\in{\bf V}$.  Next, the marginal $\mu^{\tilde f}(dy\times{\bf B})$ on $\bf V$ is the limit of the increasing (set-wise) sequence of measures 
$$\mu_n(\Gamma):=\sum_{i=1}^n \PP^{\tilde f}_{y_0}(Y_{i-1}\in\Gamma),~~\Gamma\in{\cal B}({\bf V}),~n=0,1,\ldots,$$
which satisfy the following relations:
\begin{eqnarray}
\mu_0(\Gamma)&=&0;\nonumber\\
\mu_n(\Gamma)&=&\delta_{y_0}(\Gamma)+\int_{\bf V} Q(\Gamma|y,\tilde f(y))\mu_{n-1}(dy),~~~\Gamma\in{\cal B}({\bf V}),~n=1,2,\ldots.\label{e56}
\end{eqnarray}
The proof is by induction: (\ref{e56}) is valid at $n=1$ and, if it is valid for some $n\ge 1$, then
\begin{eqnarray*}
\mu_{n+1}(\Gamma)&=&\delta_{y_0}(\Gamma)+\sum_{i=1}^n\int_{\bf V} Q(\Gamma|y,\tilde f(y))\PP^{\tilde f}_{y_0}(Y_{i-1}\in dy)\\
&=& \delta_{y_0}(\Gamma)+\int_{\bf V} Q(\Gamma|y,\tilde f(y))\mu_n(dy),~~~~~\Gamma\in{\cal B}({\bf V}).
\end{eqnarray*}
Since, according to (\ref{e55}),
$$\mu^f(\Gamma\times{\bf B})=\delta_{y_0}(\Gamma)+\int_{\bf V} Q(\Gamma|y,\tilde f(y))\mu^f(dy\times{\bf B}),~~~~~\Gamma\in{\cal B}({\bf V}),$$
$\mu_n(\Gamma)\le \mu^f(\Gamma\times{\bf B})$ for all $\Gamma\in{\cal B}({\bf V})$ and $n=0,1,\ldots$: the obvious proof is by induction w.r.t. $n$. Thus,
$$\mu^{\tilde f}(\Gamma\times{\bf B})\le \mu^f(\Gamma\times{\bf B}),~~~~~\Gamma\in{\cal B}({\bf V}).$$
We conclude that, in  turn, $f$ is a modification of $\tilde f$ on a $\mu^{\tilde f}(dx\times{\bf B})$-null subset of $\bf V$, and
$$\mu^{\tilde f}(dy\times db)=\mu^{\tilde f}(dy\times{\bf B})\delta_{\tilde f(y)}(db)=\mu^{\tilde f}(dy\times{\bf B})\delta_{f(y)}(db).$$
Repeating the above reasoning for $f$ and $\tilde f$ being swapped, we see that $\mu^{f}(\Gamma\times{\bf B})\le \mu^{\tilde f}(\Gamma\times{\bf B})$, $\Gamma\in{\cal B}({\bf V})$.

Therefore,  $\mu^{f}(dy\times{\bf B})= \mu^{\tilde f}(dy\times{\bf B})$ and hence
$$\mu^f(dy\times db)=\mu^{\tilde f}(dy\times {\bf B})\delta_{\tilde f(y)} (db)=\mu^{\tilde f}(dy\times db)$$
on ${\bf V}\times{\bf B}$ by (\ref{e55}). \hfill$\Box$
\bigskip

\par\noindent\textit{Proof of Theorem \ref{l8}.} There is a feasible and optimal in problem (\ref{PZZeqn02}) strategy $\pi^*$ in the form of a mixture of $J+1$ deterministic stationary strategies $f^1,f^2,\ldots, f^{J+1}$ according to Proposition \ref{t2}. Without loss of generality, one can say that $\pi^*\in\Pi$: as usual, just put $\pi^*_i(db|h_{i-1}):=\delta_{f^*(x_{i-1}}(db)$ if $x_{i-1}$, the last component of $h_{i-1}$ is in ${\bf V}^c$. The values of the objectives ${\cal V}(x_0,\pi^*,\bar{C}_j)$, $j=0,1,\ldots,J$ cannot increase after this modification. Moreover, we preliminarily modify the strategies $f^l$, $l=1,2,\ldots,J+1$ in the same way. Now all the measures $\mu^{f^l}(dx\times {\bf B})$ are $\sigma$-finite on ${\bf V}$ by Lemma \ref{l2}: the strategy $\pi^*$ is feasible with a finite value, hence $\int_{{\bf X}_\Delta\times{\bf B}}\bar C_j(x,b)\mu^{f^l}(dx\times db)<\infty$ for all $l=1,2,\ldots,J+1$ because $\gamma^l>0$.
The measures $\mu^{f^l}|_{{\bf V}\times{\bf B}}$ obviously belong to the set $\cal D$.

Problem (\ref{PZZeqn02})  is equivalent to the relaxed problem (\ref{SashaLp02})-(\ref{e10}) by Theorem \ref{t1}, and the restriction of the occupation measure $\mu^*:=\mu^{\pi^*}|_{{\bf V}\times{\bf B}}$ is feasible and optimal in the relaxed problem (\ref{SashaLp02})-(\ref{e10}) and also in the primal convex program (\ref{enum3}). Therefore, $L_1(\mu^*,\bar g^*)=\inf_{\mu\in{\cal D}} L_1(\mu,\bar g^*)$ according to Theorem  \ref{prop2}(b-ii).
Recall that $\bar g^*$ is the maximizer of the dual functional $h(\bar{g}).$
Since, for each $l=1,2,\ldots,J+1$, the restriction of the occupation measure $\mu^{f^l}|_{{\bf V}\times{\bf B}}$ belongs to $\cal D$, $\gamma_l>0$, $\mu^*=\sum_{l=1}^{J+1} \gamma_l\mu^{f^l}|_{{\bf V}\times{\bf B}}$, and the Lagrangian $L_1(\mu,\bar{g}^\ast)$ is linear in $\mu$ with non-negative functions $\bar C_j(\cdot,\cdot)$, $j=0,1,\ldots,J$, we see that $L_1(\mu^{f^l},\bar g^*)=\inf_{\mu\in{\cal D}} L_1(\mu,\bar g^*)$ for all $l=1,2,\ldots,J+1$.

According to Theorem \ref{t5}(b) and Theorem \ref{prop2}(a), for each $l=1,2,\ldots,J+1$,
\begin{eqnarray*}
h(\bar g^*) &=& W^*_{\bar g^*}(x_0)-\sum_{j=1}^J g^*_j d_j=\inf_{\mu\in{\cal D}} L_1(\mu,\bar g^*)=L_1(\mu^{f^l},\bar g^*)\\
&=& \int_{{\bf V}\times{\bf B}}\left(\bar C_0(x,b)+\sum_{j=1}^J g^*_j\bar C_j(x,b)\right)\mu^{f^l}(dx\times db)-\sum_{j=1}^J g^*_j d_j,
\end{eqnarray*}
so that
\begin{eqnarray*}
W^*_{\bar g^*}(x_0)&=&\int_{{\bf V}\times{\bf B}}\left(\bar C_0(x,b)+\sum_{j=1}^J g^*_j\bar C_j(x,b)\right)\mu^{f^l}(dx\times db)\\
&=&\int_{\bf V} \left(\bar C_0(x,f^l(x))+\sum_{j=1}^J g^*_j\bar C_j(x,f^l(x))\right)\mu^{f^l}(dx\times {\bf B}):
\end{eqnarray*}
the deterministic stationary strategy $f^l$ is optimal in the sense that ${\cal V}(x_0,f^l,\bar C_0+\sum_{j=1}^J g^*_j\bar C_j)=W^*_{\bar g^*}(x_0)$.
Recall that the Bellman function $W^*_{\bar g^*}(\cdot)$ is bounded by Theorem  \ref{t5}(a).
According to Lemma \ref{l6}, for each $l\in\{1,\dots,J+1\},$ $\mu^{f^l}|_{{\bf V}\times{\bf B}}=\mu^{\tilde f^l}|_{{\bf V}\times {\bf B}}$, where the deterministic stationary strategy $\tilde f^l$ is  uniformly optimal in the sense that
\begin{eqnarray*}
&&{\cal V}(x,f^l,\bar C_0+\sum_{j=1}^J g^*_j\bar C_j)\\
&=&\EE^{\tilde f^l}_x\left[\sum_{i=1}^\infty \left(\bar C_0(X_{i-1},B_i)+\sum_{j=1}^J g^*_j\bar C_j(X_{i-1},B_i)\right)\right]\\
&=&\inf_\pi \EE^{\pi}_x\left[\sum_{i=1}^\infty \left(\bar C_0(X_{i-1},B_i)+\sum_{j=1}^J g^*_j\bar C_j(X_{i-1},B_i)\right)\right]= W^*_{\bar g^*}(x),~~~x\in{\bf V}.
\end{eqnarray*}
Now $\tilde f^l\in{\bf F}$ by Proposition 9.12 of \cite{Bertsekas:1978}, and $\pi^*$ is a mixture of $\tilde f^1,\tilde f^2,\ldots,\tilde f^{J+1}$. $\hfill\Box$
 
 \subsection*{Acknowledgment}
Part of this work was carried out during the visit of the second author to Sun Yat-sen University, China, and thus, this work is partially supported by the National Natural Science Foundation of China (No.72471252).

\end{document}